\documentclass{article}
\usepackage[top=1in, bottom=1in, left=1in, right=1in, a4paper]{geometry}

\usepackage{hyperref}
\usepackage{color}
\usepackage[x11names,dvipsnames,table]{xcolor}
\usepackage{colortbl}
\usepackage{graphicx}
\usepackage{cite}
\usepackage{array}
\usepackage{booktabs}
\usepackage{listings}
\usepackage{longtable}
\usepackage{amssymb}

\usepackage{authblk}

\begin{document}

\title{Reinforcement learning for graph theory, \\ I. Reimplementation of Wagner's approach\footnote{%
       This work was supported and funded by Kuwait University Research Grant No. SM05/22.}}

\author[a]{Mohammad Ghebleh\thanks{Corresponding author.}}
\author[a]{Salem Al-Yakoob}
\author[a]{Ali Kanso}
\author[b]{Dragan Stevanovi\'c\thanks{On leave from Mathematical Institute of Serbian Academy of Sciences and Arts.}}
\affil[a]{Department of Mathematics, Faculty of Science, Kuwait University, Safat, Kuwait}
\affil[b]{College of Integrative Studies, Abdullah Al Salem University, Khaldiya, Kuwait}
\affil[ ]{}
\affil[ ]{mghebleh@gmail.com, salem.alyakoob@ku.edu.kw, akanso@hotmail.com, dragan.stevanovic@aasu.edu.kw}

\date{}

\maketitle

\begin{abstract}
We reimplement here the recent approach of Adam Zsolt Wagner [arXiv:2104.14516],
which applies reinforcement learning to construct (counter)examples in graph theory,
in order to make it more readable, more stable and much faster.
The presented concepts are illustrated by constructing counterexamples 
for a number of published conjectured bounds for the Laplacian spectral radius of graphs.
\end{abstract}



\section{Introduction}

We are witnessing yet another computational revolution 
through proliferation of applications of artificial intelligence (AI) methods in everyday tasks. 
Besides being very publicly scrutinised 
through the ongoing development of self-driving vehicles, 
defeating world champions in Go, 
achieving superhuman performance in recognising images and playing video games, 
or generating multimedia content, 
AI models are trained and used to an even greater extent in science, engineering and technology 
to provide surrogate models for predicting values of functions that are hard to compute or simulate numerically. 
Our goal here is to elaborate on a recent application of a particular AI method—--reinforcement learning—--%
in construction of (counter)examples in graph theory.

A few specialised software packages have been used 
as auxiliary tools in graph theoretical research for four decades already, 
either to help with posing new conjectures 
or to help refute existing conjectures by providing counterexamples. 
Two earliest examples of such packages are GRAPH and Graffiti, 
which were special types of expert systems. 
GRAPH, written by Dragoš Cvetković and Laszlo Kraus in the 1980s \cite{1,2,3}, 
provided a closed environment for visually editing individual graphs and computing their invariants, 
and implemented certain AI methods for automatically proving simplest theorems in graph theory. 
Graffiti, written by Siemion Fajtlowicz in 1986 \cite{4,5,6,7}, 
was geared more toward automatic conjecture making than enabling researcher to test his own conjectures. 
Graffiti was used to produce many conjectures, 
some of which attracted attention of well-known graph theorists \cite{8}. 
It is interesting to note that these old packages
have been rewritten and modernised a few times in the forms of 
newGRAPH \cite{9,10,11}, Graffiti.pc \cite{12,13,14,15} and Grinvin \cite{16,17}. 
A relatively recent addition to this group is graph6java \cite{18,19}, 
which builds upon the power of existing programs for exhaustive generation of different types of graphs 
and is capable, for a selected graph invariant, to answer the questions 
such as which graphs have extremal values or which pairs of graphs have equal values of this invariant. 
Another approach is taken by AutoGraphiX (AGX), 
written by Gilles Caporossi and Pierre Hansen in the 2000s \cite{20,21,22,23,24}. 
AGX considered expressions made up of graph invariants 
as instances of an optimisation problem and 
consequently applied variable neighbourhood search (VNS) metaheuristic 
to search for extremal graphs for these expressions 
over a search space that consists of a predefined set of graphs 
(which is too big to perform an exhaustive search in it). 
AGX has been used to refute some of the conjectures of Graffiti, and 
in a lengthy series of papers its authors posed several new conjectures as well.

\begin{figure}[ht]
\begin{center}
\includegraphics[width=0.75\textwidth]{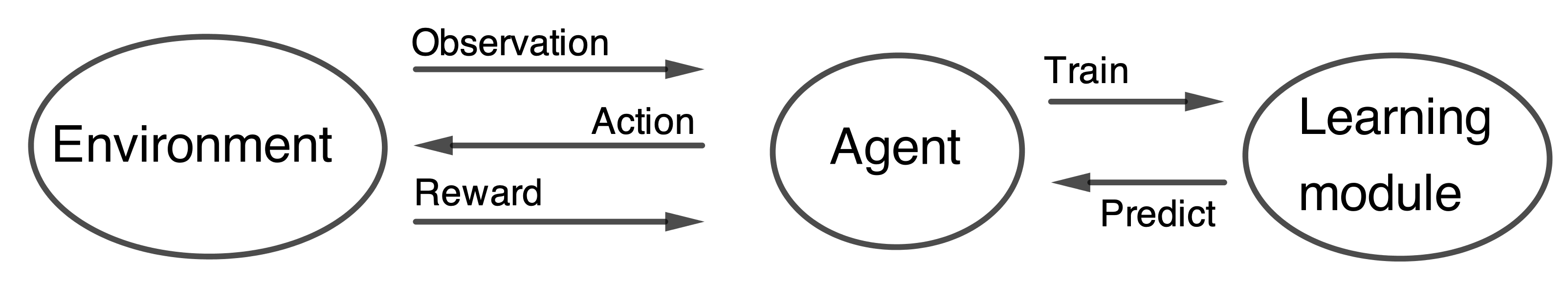}
\end{center}
\caption{General process of reinforcement learning.}
\label{fig-rl-principle}
\end{figure}

Most recently, in the pioneering work~\cite{25}
Adam Zsolt Wagner used reinforcement learning (RL) 
to construct counterexamples for several conjectures in graph theory.
Although RL is present for a long time in AI community,
it experienced recent rise in popularity after it was used by DeepMind in programs AlphaGo, AlphaZero and MuZero
that managed to reach superhuman level of performance in go, chess and visually complex Atari games, 
all without knowing the game rules in advance \cite{26}.
RL is based on an interplay between an agent and an environment (see Fig.~\ref{fig-rl-principle}):
at each step,
the agent receives the observation on the state of the environment,
performs an action in the environment,
for which it then receives the reward from the environment.
The goal of the RL agent is to learn how to maximise the cumulative reward received from the environment.

In Wagner's approach, 
simple graphs are constructed by an agent as a sequence of 0-1 actions
which represent the entries of the upper part of the adjacency matrix.
The environment responds with zero intermediate rewards,
while the final reward, when the graph is fully constructed, is equal to a specified graph invariant.
Setting the final reward to represent the difference of two sides of a conjectured inequality between graph invariants, 
Wagner \cite{25} was able to find counterexamples for several published conjectures. 
Among others, this includes
a 19-vertex counterexample for an AGX conjecture on the sum of the matching number and the spectral radius~\cite{20},
for which one of the authors previously showed the existence of counterexamples on more than 600 vertices~\cite{27}.
One of the authors~\cite{28} also tested Wagner's approach
on a recent conjecture of Akbari, Alazemi and Andjeli\'c~\cite{29}.
Although RL did not manage to explicitly find counterexamples in this case, 
it was still very useful by suggesting their proper structure.

Our goal here is to offer the reimplementation of Wagner's approach
that is more readable, more stable, much faster and as a result, more useful.
Wagner's initial RL implementation runs only for very specific combinations of python and tensorflow versions,
which is probably a consequence of pre-existing code that it is based on (see \cite{30}).
To avoid this, we reimplement his approach here from scratch.
We separate the learning agent from the reward computation,
so that interested researchers need only define the final reward in a separate python file 
in order to adapt and reuse this implementation.
We pay particular attention to computation performance,
since RL must compute invariants for at least several hundred thousands of graphs
before reaching a satisfactory level of convergence in its learning.
While final rewards may certainly be computed by using networkX and/or numpy,
our experience shows that computing them by calling Java code directly from Python offers a significant speedup.
For example, 
computation of eigenvalues of graphs on 20--30 vertices 
with graph6java (which internally uses EJML library \cite{31}) is
3--5 times faster than computing them with numpy.
Overall, the new implementation enables interested researchers to see learning results 
in a matter of minutes instead of hours or even days.

In the next section we briefly overview the principles used in Wagner's approach
and the details of our reimplementation, which also serves as a short manual for its use.
In Section~\ref{sc-autolaplabounds} we then apply it
to find counterexamples to a number of older conjectured bounds on the Laplacian spectral radius of graphs from~\cite{BrHS}.

\section{Wagner's simple graph environment and the cross-entropy method}

Wagner's most profound observation from~\cite{25} is that 
construction of a simple graph and computation of its invariant
can be treated as a sequence of consecutive observations, actions and rewards.
For an $n$-vertex graph~$G$,
let $L\colon l_1,\dots,l_{{n\choose 2}}$ be the part of adjacency matrix above the diagonal,
listed in the row-wise order.
Since $G$ is fully determined by~$L$,
the agent constructs $G$ simply by issuing actions equal to the consecutive elements of~$L$.
All rewards given by the environment before $G$ is fully constructed are equal to zero,
while the final reward is equal to the graph invariant~$r(G)$ that is of interest to us.

\begin{figure}[ht]
\begin{center}
\begin{minipage}{0.2\textwidth}
\begin{center}
\includegraphics[width=\textwidth]{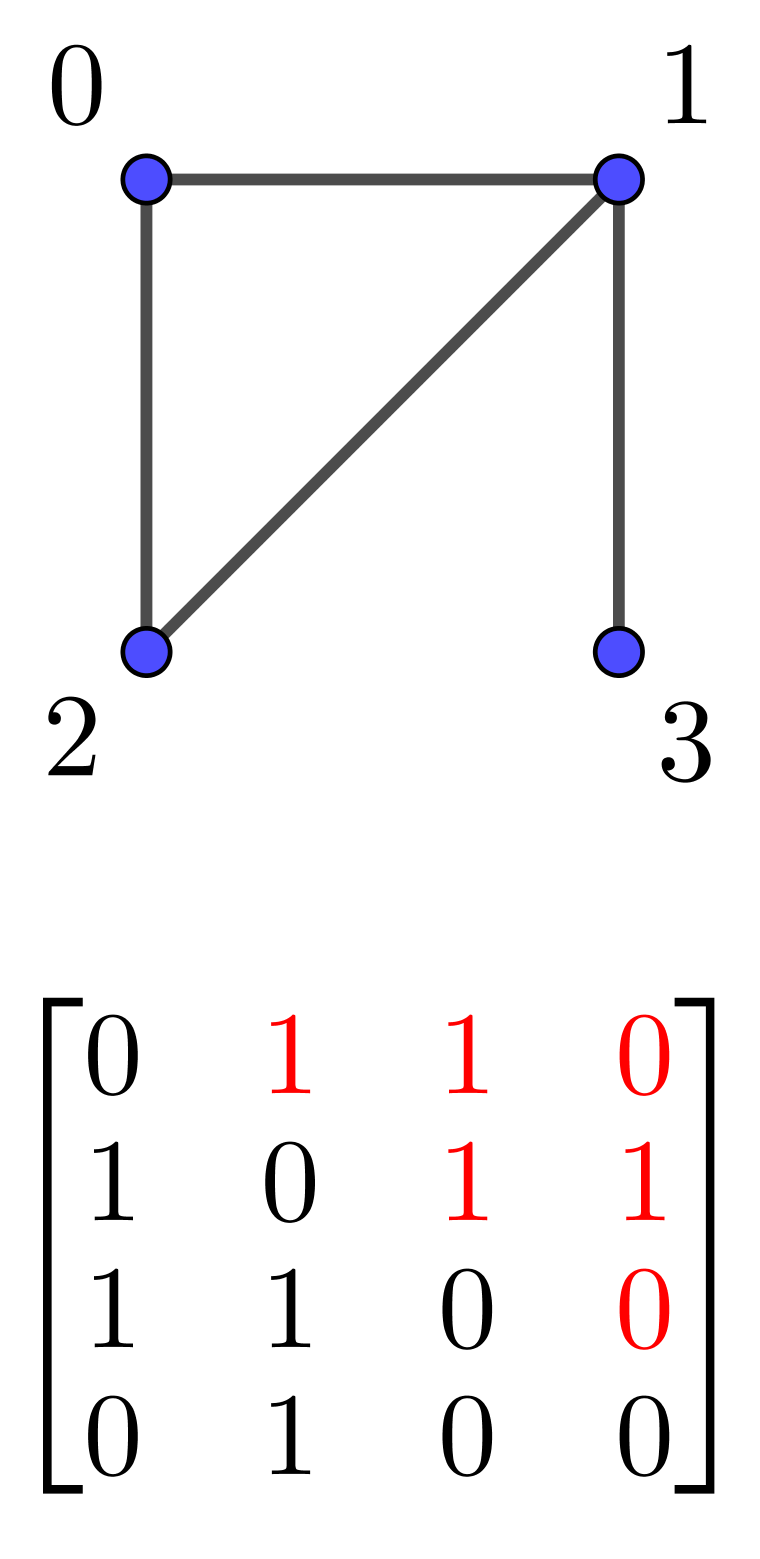}
\end{center}
\end{minipage}
\qquad\qquad
\begin{minipage}{0.5\textwidth}
\begin{center}
\begin{tabular}{ccc}
\toprule
Environment & Agent & Environment \\
\midrule 
Observation & Action & Reward \\
\midrule 
\textcolor{red}{0}00000 \textcolor{red}{1}00000 &   & \\
                                                & 1 \\
1\textcolor{red}{0}0000 0\textcolor{red}{1}0000 &   & 0 \\
                                                & 1 \\ 
11\textcolor{red}{0}000 00\textcolor{red}{1}000 &   & 0 \\
                                                & 0 \\
110\textcolor{red}{0}00 000\textcolor{red}{1}00 &   & 0 \\
                                                & 1 \\
1101\textcolor{red}{0}0 0000\textcolor{red}{1}0 &   & 0 \\
                                                & 1 \\
11011\textcolor{red}{0} 00000\textcolor{red}{1} &   & 0 \\
                                                & 0 \\
110110                  000000                  &   & $r(G)$ \\
\bottomrule
\end{tabular}
\end{center}
\end{minipage}
\end{center}
\caption{Illustration of the sequence of observations, actions and rewards issued 
         while constructing the graph shown in top left part.
         Entries of the adjacency matrix above the diagonal, which yield the sequence~$L$, are shown in red.
         Red-coloured entries of observations represent 
         the entry in the left part which will store the next action and
         the single 1 in the right part which indicates the index of this entry.
         Note that both the list of actions issued and the left part of the last observation are equal to~$L$.}
\label{fig-graph-construction}
\end{figure}

On the other hand, environment observations have to support the goal of the agent
to learn which actions should be issued at which stage of graph construction
in order to maximise the final (=cumulative) reward.
Hence each observation has to inform the agent what part of the graph has been already constructed so far.
For this reason,
each observation represents the current state of~$L$,
followed by the one-hot encoding of the next entry of~$L$ to be determined.
The initial observation, before any action is issued, is thus 
$$
\underbrace{0\dots0}_{{n\choose 2}}\,\underbrace{10\dots0}_{{n\choose 2}}.
$$
Subsequent actions are collected in the first part of the observation which builds up to~$L$,
while the single~$1$ in the second part travels from left to right,
indicating the index of the entry of~$L$ that will store the next action,
until it should pass over the right end, 
which indicates that the graph is fully constructed and that $r(G)$ can be computed.
This process of interaction between the agent and the environment
while constructing a small graph is illustrated in Fig.~\ref{fig-graph-construction}.

\begin{figure}[ht]
\begin{center}
\includegraphics[width=0.95\textwidth]{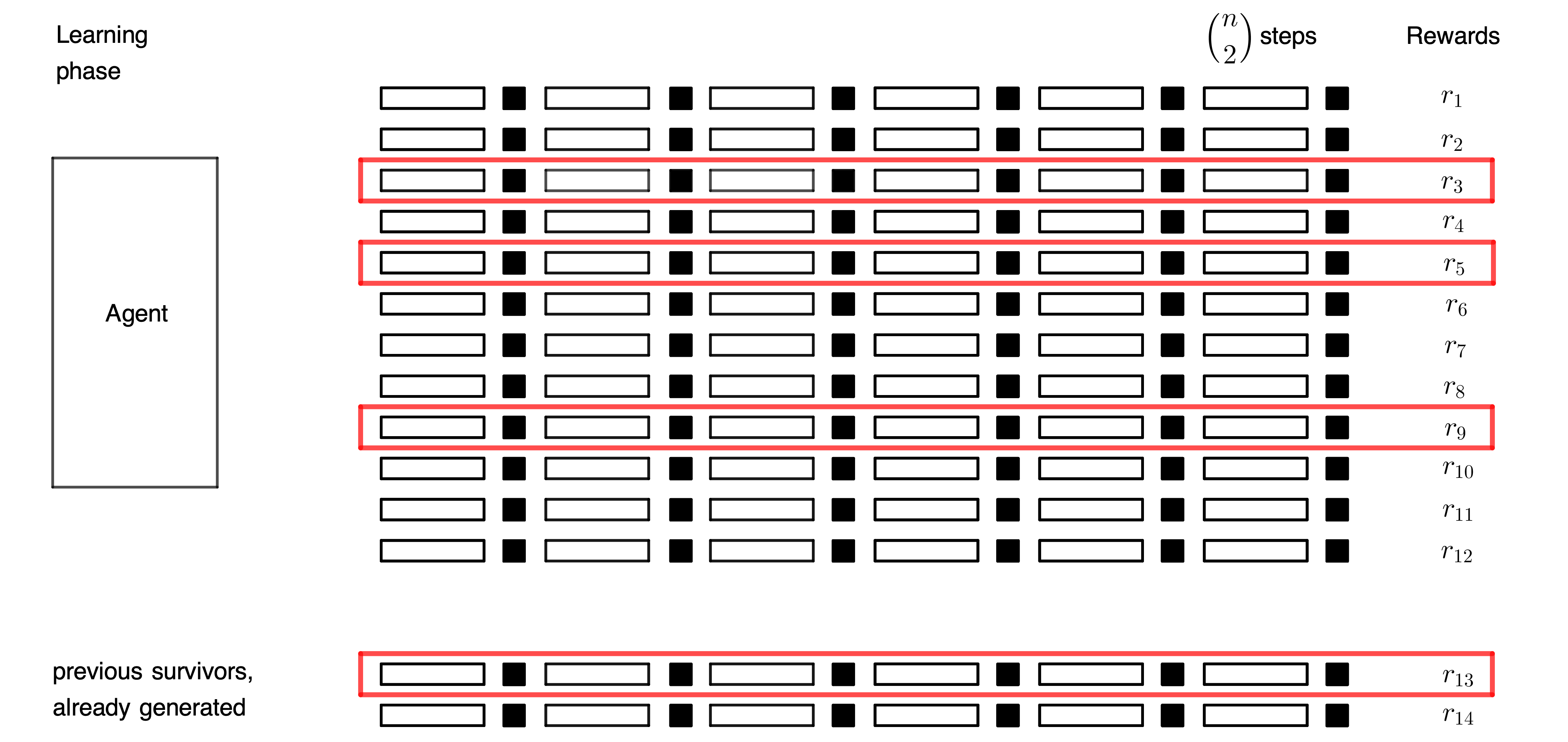}
\end{center}
\caption{Illustration of the learning phase in the cross-entropy method.
White-coloured rectangles symbolise observations of graph environments,
while black-coloured boxes represent subsequent actions issued by the agent.
Construction of each graph is completed after ${n\choose 2}$ actions,
after which the rewards $r_i$ are computed.
Constructed graphs with a specified percentage of largest rewards are selected (framed in red rectangles).
The agent's neural network is retrained on the set of observation/actions pairs used in construction of most rewarding graphs,
with the aim of minimising the cross-entropy 
between the action probability distributions that the neural network provides for these observations 
and the actual actions used in their construction.}
\label{fig-cema-learning}
\end{figure}

Wagner~\cite{25} used the \emph{cross-entropy method}~\cite{cema1,cema2} 
to enable the agent to learn how to construct graphs with large rewards.
This is one of the simplest reinforcement learning methods,
which usually works well for environments 
that do not require discovery of complex, multistep policies~\cite[Chapter~4]{lapa}.
In this method the agent employs as a learning module
the neural network (initialised randomly at the beginning)
which accepts a graph observation as the input
and outputs a probability distribution over the set of possible actions (0 and 1 in this case),
according to which the agent selects the next action.
The cross-entropy method successively iterates through the \emph{generation} and the \emph{learning} phases.
In the generation phase,
the agent uses the neural network to simultaneously construct a batch of $N$~graphs, for a fixed value of~$N$.
While each graph in the batch starts with the same initial observation $\underline{0\dots0}\,\underline{10\dots0}$,
the agent's actions are selected at random according to the probability distribution 
provided by the neural network for each observation,
so that the fully constructed graphs in the batch will not necessarily be mutually isomorphic.
In the learning phase,
the agent selects $L$~graphs (for some $L<\!\!<N$) with the largest rewards and
trains its neural network on the pairs of observations and actions used in their construction
(see Fig.~\ref{fig-cema-learning}).
The goal of training is to minimise the cross-entropy between the actions used and
the probability distributions that the network outputs for the corresponding observations.
After training is done,
the agent selects $S$~graphs (for some $S<L$) with the largest rewards
that are forwarded as survivors to the next generation phase.

Note that the only purpose of rewards in the cross-entropy method is to select the subsets of elite graphs,
as the actual value of rewards is not used in training of the agent's underlying neural network.
Nevertheless, the network learns the common characteristics of the structure of elite graphs,
so that the graphs constructed in the following generation phases will tend to have higher rewards.
Survivor graphs serve to speed up the learning process in the beginning,
directing neural network toward imitating construction of particularly rewarding graphs.
However, as the learning process converges,
the graphs constructed in each new generation will begin to share more and more similarities,
to the point when majority of them will actually be mutually isomorphic to the same locally optimal graph.
At such moment either the learning should be stopped or 
the actions issued by the agent should be additionally randomised in an effort to avoid this local optimum.

\subsection{New reimplementation}

\lstset{language=Python}

Our reimplementation of the cross-entropy method employing Wagner's simple graph environment 
is freely available online at
\begin{center}
\url{github.com/dragance106/cema-for-graphs}
\end{center}
and it consists of three independent files: 
\begin{itemize}
\item \verb|cema_train_simple_graph.py|,
\item \verb|graph6java.jar| and
\item \verb|training_runner.py|,
\end{itemize}
that we briefly discuss here.

The file \verb|cema_train_simple_graph.py| contains implementation of
the methods for successively generating batches of simple graphs and 
training the agent's neural network to maximize the reward function.
The main method in this file is:
\begin{lstlisting}[frame=tb,language=Python]
train(compute_reward,
      n=20,
      batch_size=200,
      num_generations=1000,
      percent_learn=90,
      percent_survive=97.5,
      neurons=[72,12],
      learning_rate=0.003,
      act_rndness_init=0.005,
      act_rndness_wait=10,
      act_rndness_mult=1.1,
      act_rndness_max=0.025,
      verbose=True,
      output_best_graph_rate=25)
\end{lstlisting}
Let us quickly discuss the parameters of this method:
\begin{itemize}
\item compute\_reward:
      denotes the method that computes a real-valued reward for each constructed graph.
      This method has to be provided by the user and 
      it has the signature \verb|my_reward(n, A)|,
      where \verb|n| is the number of vertices and \verb|A| is the adjacency matrix;
\item n: 
      the number of vertices of graphs to be constructed;
\item batch\_size: 
      the number of graphs in each batch;
\item num\_generation: 
      the number of generations over which to train the agent;
\item percent\_learn: 
      the best (100-percent\_learn) percents of graphs from each batch are used to train the agent's neural network;
\item percent\_survive: 
      the best (100-percent\_survive) percents of graph from each batch are transferred to the next generation as survivors;
\item neurons: 
      the list determining the numbers of neurons in hidden layers of the agent's neural network.
      The neural network has $n(n-1)$ inputs to accommodate a simple graph observation,
      and two outputs representing raw scores for the possible actions (0=skip this edge, 1=add this edge).
      The structure of hidden layers (together with their number) is defined by this parameter;
\item learning\_rate:
      the learning rate of the agent's neural network;
\item act\_rndness\_init: 
      the initial value for the action randomness act\_rndness.
      In order to increase exploration of different graphs and
      avoid being stuck in a local optimum from overexploitation of acquired knowledge,
      the agent will issue random actions at this rate,
      i.e., this percent of adjacency matrix entries will be random;
\item act\_rndness\_wait:   
      the number of generations without an increase in the maximum reward
      to wait before act\_rndness should be increased;
\item act\_rndness\_mult:   
      the factor used to multiply act\_rndness 
      when there are no increases in the maximum reward for act\_rndness\_wait generations;
\item act\_rndness\_max:    
      the maximum allowed act\_rndness value,
      since we do not want to have too many random edges in constructed graphs;
\item verbose:            
      True or False, describing whether to print on the console the summary information for each generation:
      generation number, maximum reward seen, 
      reward used to select survivors, 
      reward used to select graphs for learning,
      processing time for the generation, and act\_rndness value
\item output\_best\_graph\_rate:
      the number of generations at which to produce the drawing of the best graph so far,
      which is then reported in the external runs/event file for TensorBoard consumption.
\end{itemize}

When the learning starts to converge,
the train method tries to prevent mutual isomorphism of constructed graphs
by randomising a certain percentage of actions issued by the agent.
This percentage is represented by the internal variable act\_rndness,
while the four act\_rndness\_\dots parameters above serve to control it
in a way that is inspired by the variable neighborhood search~\cite{VNS}:
whenever act\_rndness\_wait generations pass without improving the maximum reward,
act\_rndness is multiplied by act\_rndness\_mult (with the waiting counter reset),
but it is not allowed to surpass the value of act\_rndness\_max.
On the other hand, as soon as the maximum reward increases in the next generation,
act\_rndness is returned back to act\_rndness\_init.
Such varying action randomness can somewhat postpone convergence towards generating isomorphic graphs,
depending on the combination of these four parameters.

After training is completed, 
the train method returns the maximum reward obtained,
as well as the adjacency matrix of a graph that attains it.
However, in order to check up on the progress of learning during training,
the train method reports both short textual summary information on the console after each generation 
(when verbose is set to True) and
writes this data, together with a drawing of the best graph each output\_best\_graph\_rate generations,
to an external event file in the runs subfolder. 
This event file may be visualised already while it is being populated during training,
by starting TensorBoard application from the same folder with:
\begin{verbatim}
tensorboard --logdir runs
\end{verbatim}
and then opening the web page indicated by TensorBoard in the browser.
In this way, 
the user can effectively visualise the evolution of both the maximum reward and the best graphs during training,
which helps her/him to judge whether training is going in a good direction or whether it should be stopped or restarted.
Note, however, that there is not a simple recipe for determining the proper combination of training parameters---%
at the end, it boils down to experimentation through trial and error.

The values indicated after the equality signs in the signature of the train method above
represent the default values for the corresponding parameters
(which need not be mentioned in the method call if they are not changed).
Hence the minimum way of calling the train method is as follows:
\begin{lstlisting}[language=Python]
from cema_train_simple_graph import train
r, A = train(compute_reward=my_reward)
\end{lstlisting}
where my\_reward is a method that has to be provided by the user.
This method accepts the number of vertices \verb|n| and the adjacency matrix \verb|A| for each constructed graph,
and returns a real number that represents its reward.
If we are looking for a counterexample to a conjectured inequality $a(G)\leq b(G)$,
Wagner's suggestion~\cite{25} was to return the reward $a(G)-b(G)$,
since in this way we can easily recognise counterexamples through the positive value of their reward.

While it is quite standard in the Python community
to use networkX~\cite{networkx} and numpy~\cite{numpy} for computation of graph invariants,
we have found that using our existing Java framework graph6java~\cite{18,19} directly from Python
offers significant performance speedup---between three and five times---without the need for any compilation.
This is enabled by the package JPype~\cite{33},
which starts the Java virtual machine parallel to the Python virtual machine and sharing the same memory,
so that with a small overhead, one can call Java methods directly from the Python code.
In order not to interfere with the implementation details in \verb|cema_train_simple_graph.py|,
it is customary to write the reward method and start the training from another python file.
Here is a minimal example that
implements the reward for the conjectured inequality
$$
\mu\leq\max_{v\in V(G)} \sqrt{\frac{4d_v^3}{m_v}}
$$
from~\cite{BrHS},
where $\mu$ is the largest Laplacian eigenvalue of the graph,
$d_v$ is the degree of vertex~$v$ and $m_v$ is the average degree of the neighbours of~$v$:

\begin{lstlisting}[frame=tb, language=Python]
import jpype.imports
from jpype.types import *
jpype.startJVM(classpath=['*'], convertStrings=False)
from graph6java import Graph

import numpy as np
import math
MINUS_INF = -1000000     # reward signifying unwanted graphs

def auto_lapla_1(n, A):
    g = Graph(JInt[:,:](A))
    if g.numberComponents()>1:
        return MINUS_INF

    mu = max(g.Lspectrum())
    deg = np.array(g.degrees())
    avd = np.array(g.averageDegrees())

    return mu - max(np.sqrt(4*deg*deg*deg/avd))

from cema_train_simple_graph import train
r, A = train(compute_reward=auto_lapla_1)

jpype.shutdownJVM()
\end{lstlisting}

The above code makes an important assumption that Java archive \verb|graph6java.jar|
(available at \url{github.com/dragance106/cema-for-graphs}) is located in the current folder.
The list of invariants currently implemented in its Graph class is given in Appendix~\ref{sc-graph6java}.
When needed, new invariants can be added to the Graph class directly in the graph6java source files,
available at \url{github.com/dragance106/graph6java},
and the newly compiled jar file can then be used instead.

Finally, the file \verb|training_runner.py|
contains many further examples of implemented reward functions and
different ways of inviting the train method.

\section{Automated conjectures on the Laplacian spectral radius of graphs}
\label{sc-autolaplabounds}

Brankov, Hansen and Stevanovi\'c observed in~\cite{BrHS} that 
a large number of upper bounds on the Laplacian spectral radius~$\mu$ of graphs published in the literature so far
have very similar form: 
they are a maximum taken either over the vertices or over the edges of 
functions that depend on the degree~$d_v$ of the vertex and the average degree $m_v$ of the neighbors of~$v$,
and that evaluate to~$2x$ when all the $d_v$ and $m_v$ terms in the functional expression are replaced by~$x$.
Based on this observation,
they proposed a procedure for automatically creating new conjectured bounds of increasing complexity.
Out of 361 vertex-maximum bounds and 1138 edge-maximum bounds of small complexity generated in~\cite{BrHS},
it turned out that 190 vertex-maximum and 297 edge-maximum bounds hold
for connected graphs with up to 9 vertices, as well as on stars and windmills (triangles sharing one common vertex).
This is a fairly large ratio of the total number of generated bounds,
suggesting that the way of generating the bounds is likely meaningful.
The 68 conjectured bounds from~\cite{BrHS}, having either small complexity or being most interesting,
are listed here in Appendix~\ref{sc-autolaplabounds-list}.


Here we apply the previously described reimplementation of 
Wagner's approach to reinforcement learning on graphs
to all 68 conjectured automated bounds on the Laplacian spectral radius of graphs.
The training was run mostly for graphs on 20~vertices, ocassionally going up to 24~vertices,
and it managed to find counterexamples for 25 of the conjectured bounds:
for 8 out of 32 vertex-maximum bounds and for 17 out of 36 edge-maximum bounds,
as indicated in Table~\ref{tb-counterexamples}.
Tensorboard's event files, from which one can see the actual counterexamples
and the evolution of maximum, learning and surviving rewards for all of these conjectures,
can be downloaded from~\cite{zenodo}.

\begin{figure}
\centering
\includegraphics[width=0.45\textwidth]{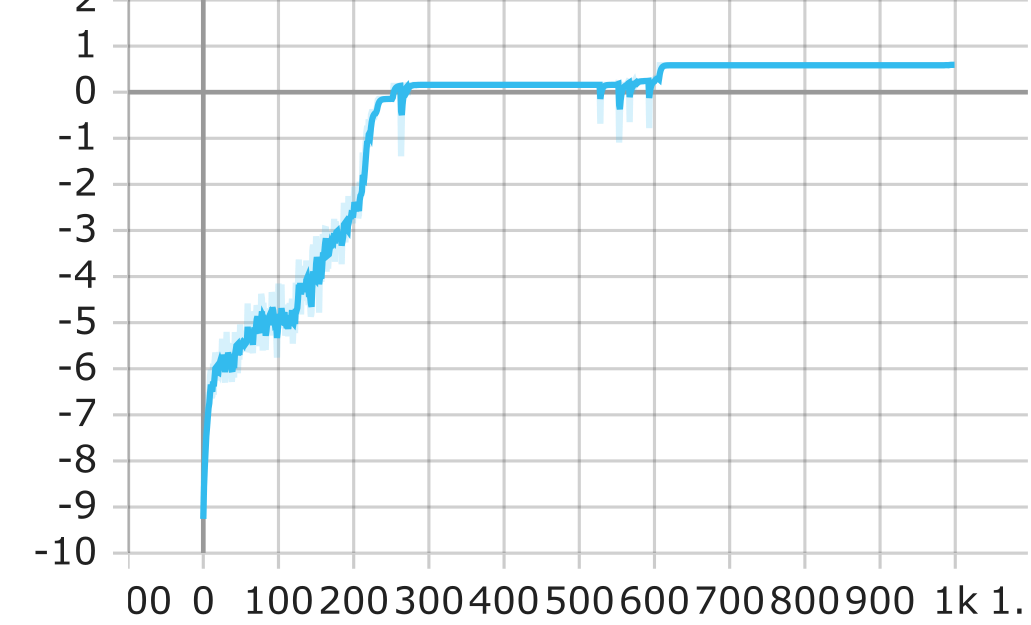}
\qquad
\includegraphics[width=0.45\textwidth]{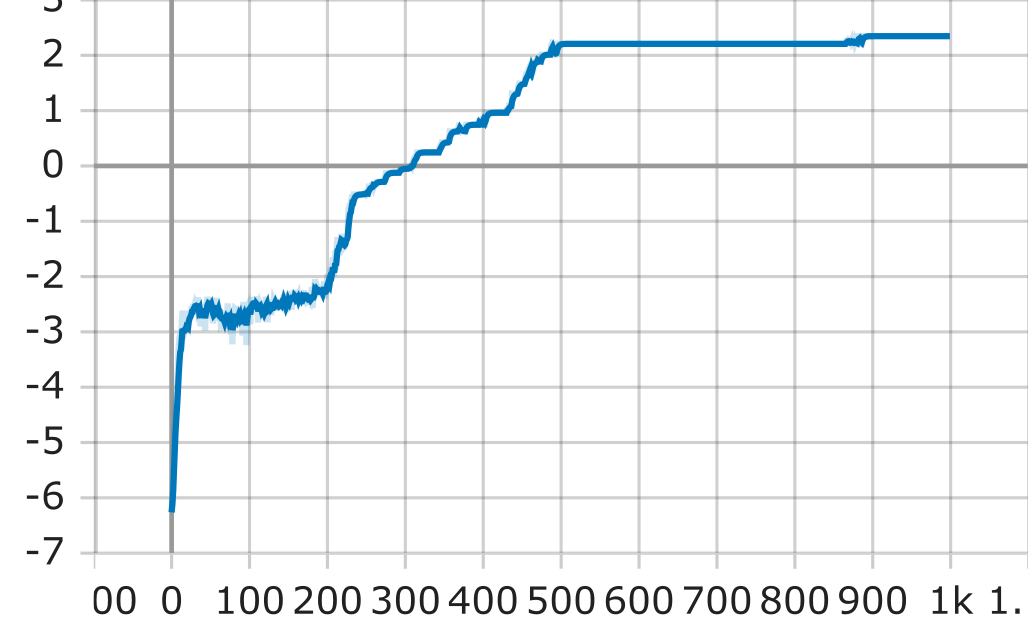}
\caption{Evolution of learning rewards for the conjecture 
         {\bf (left)} 31: $\mu\leq\max_{v\in V} \frac{4m_v^2}{m_v+d_v}$
         and {\bf (right)} 65: $\mu\leq\max_{v_i\sim v_j} \frac{(m_i+m_j)(d_im_i+d_jm_j)}{2m_im_j}$.}
\label{fig-learning-evolution}
\end{figure}

\begin{figure}
\centering
\begin{tabular}{ccc}
\includegraphics[width=0.225\textwidth]{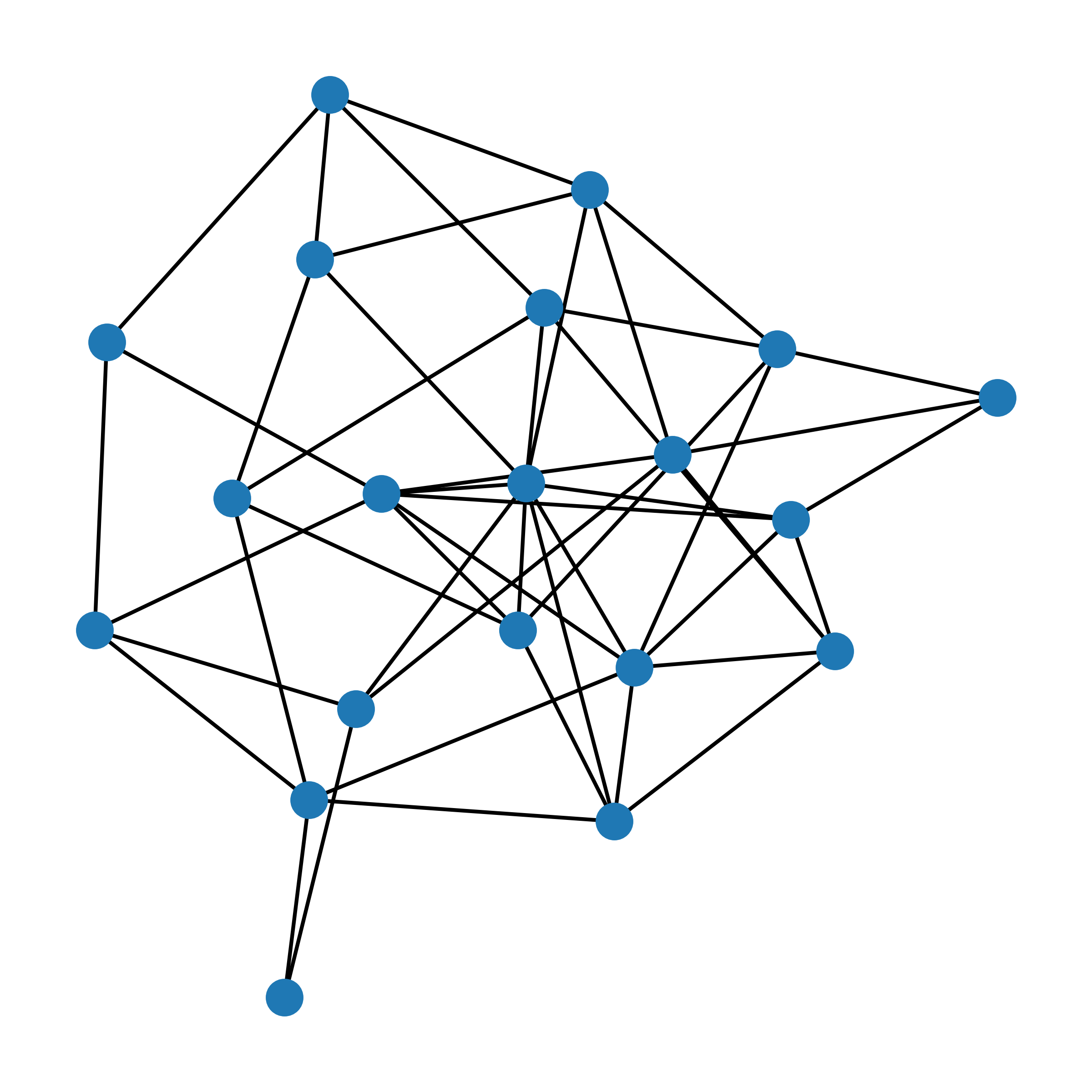} 
  & \includegraphics[width=0.225\textwidth]{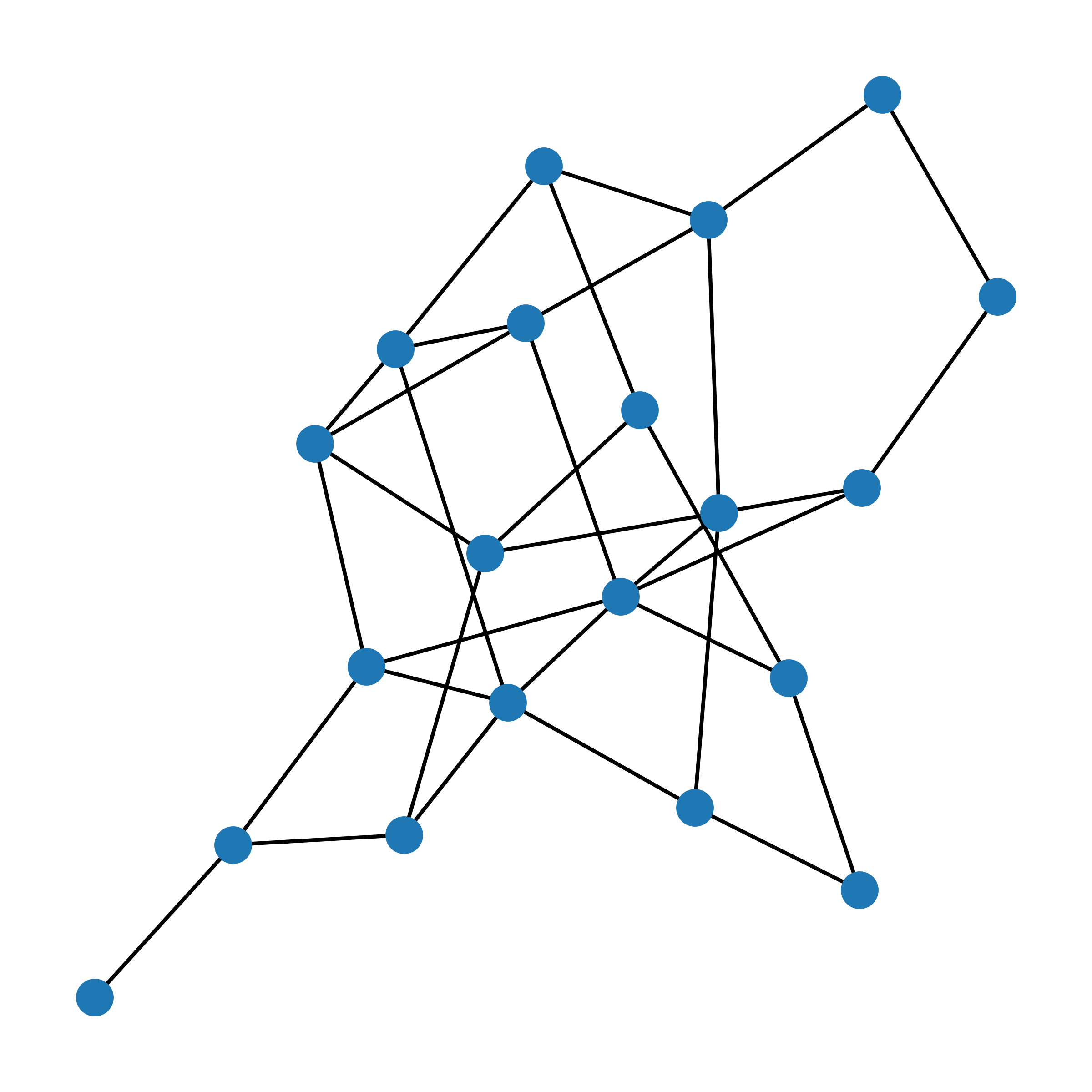} 
  & \includegraphics[width=0.225\textwidth]{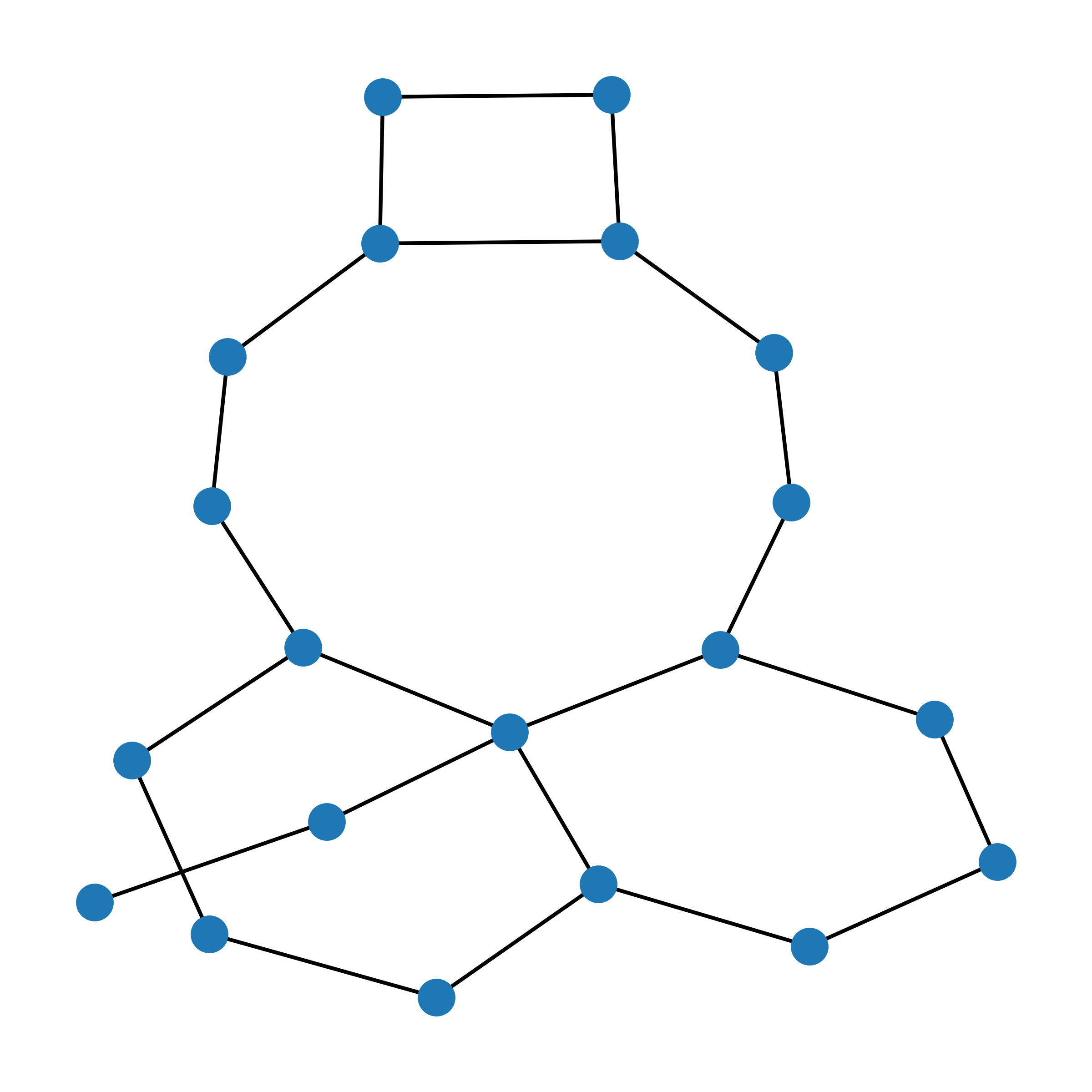} \\
gen=50, reward=-3.106
  & gen=75, reward=-2.371
  & gen=250, reward=0.1517 \\
\includegraphics[width=0.225\textwidth]{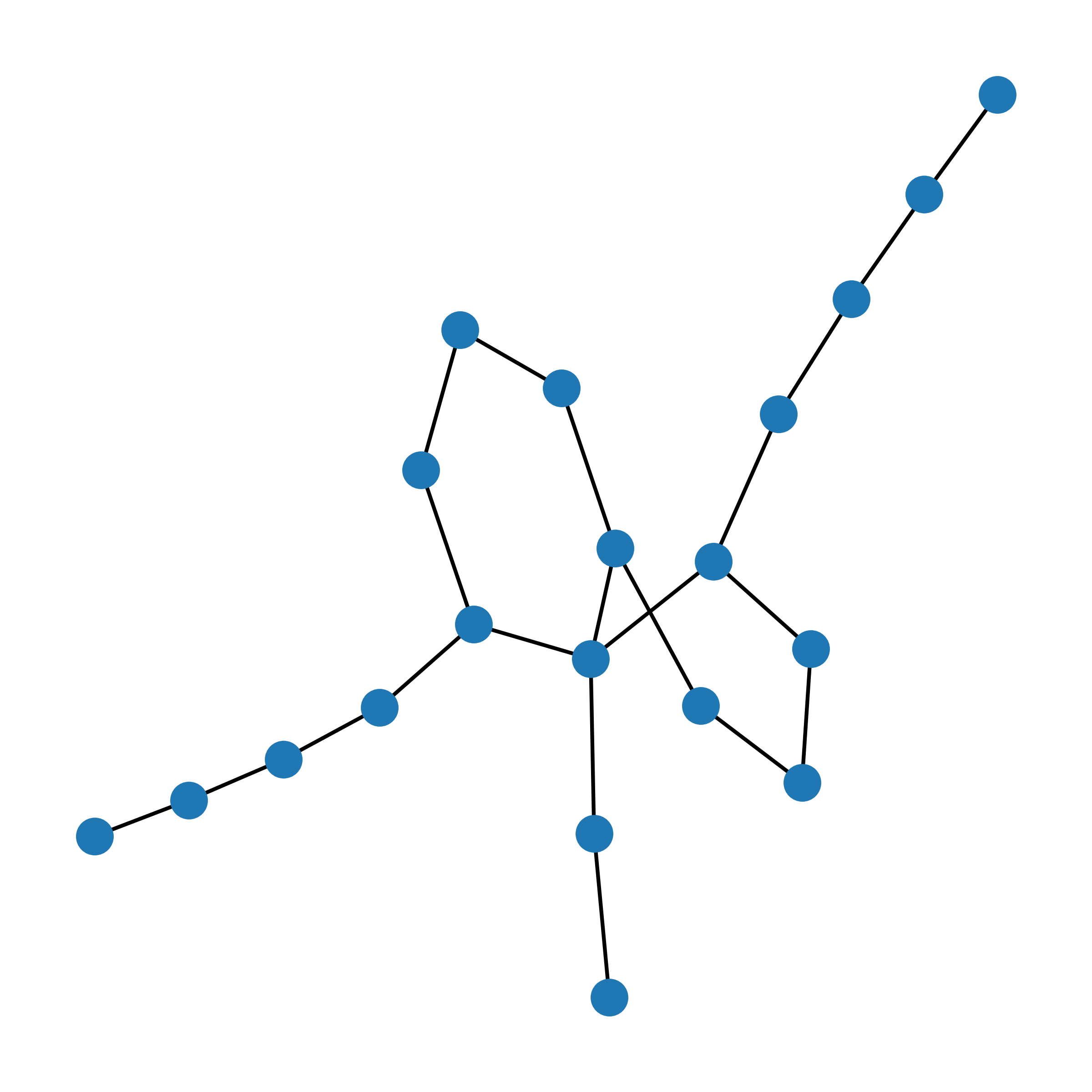} 
  & \includegraphics[width=0.225\textwidth]{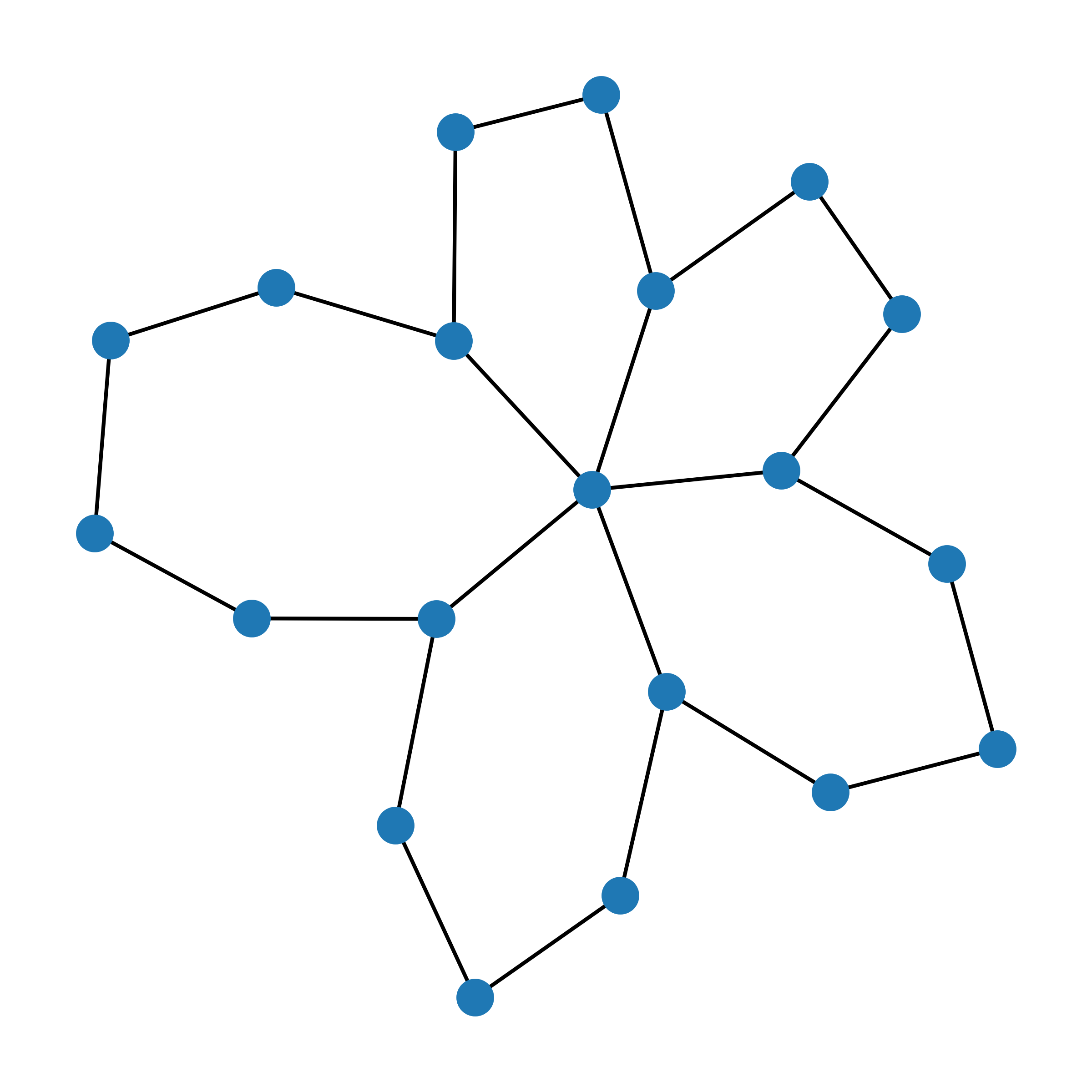} \\
gen=500, reward=0.1588
  & gen=600, reward=0.5840
\end{tabular}
\caption{Evolution of graphs with maximum rewards for the conjecture 
         31: $\mu\leq\max_{v\in V} \frac{4m_v^2}{m_v+d_v}$.}
\label{fig-evolution-31}
\end{figure}

\begin{figure}
\centering
\begin{tabular}{ccc}
\includegraphics[width=0.225\textwidth]{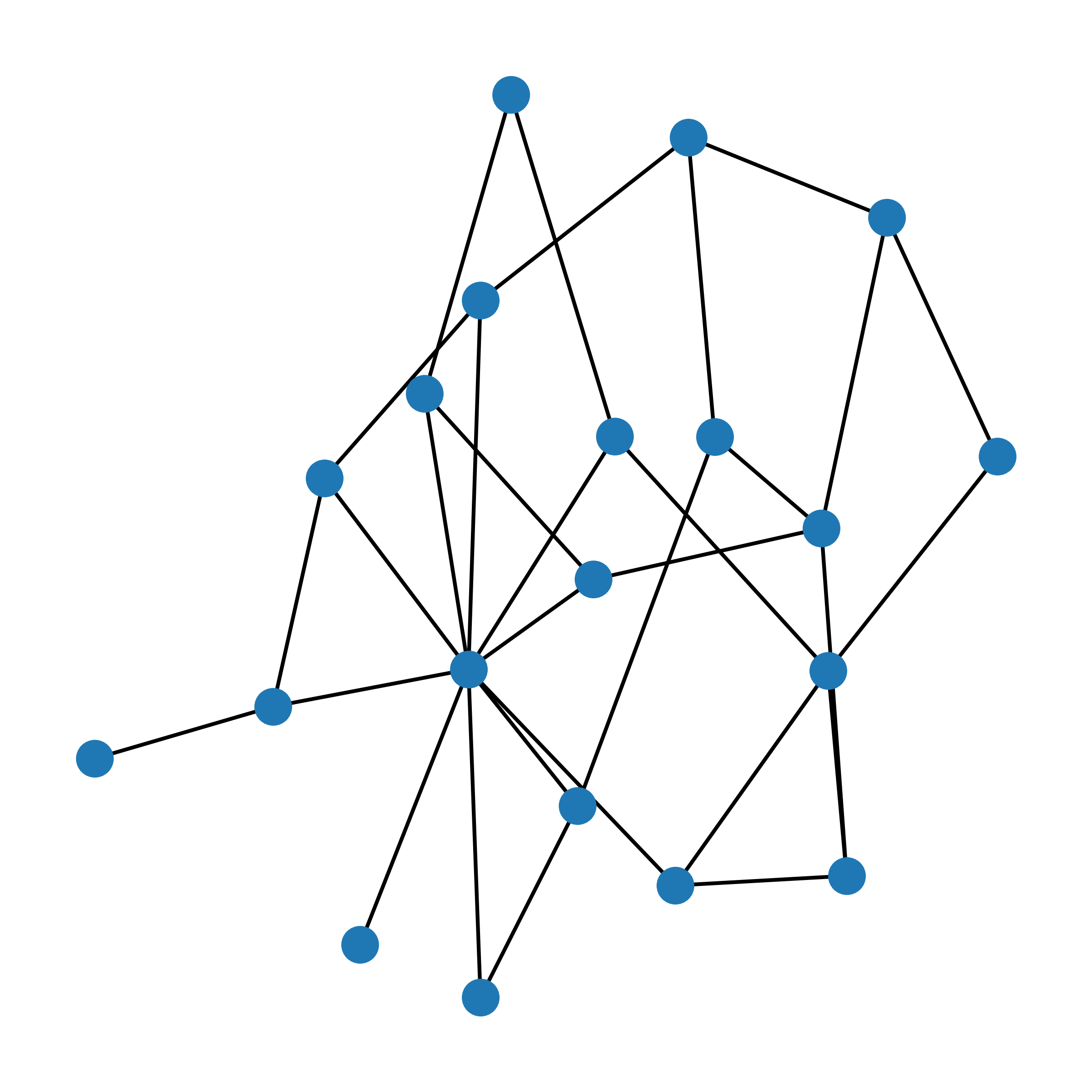}
  & \includegraphics[width=0.225\textwidth]{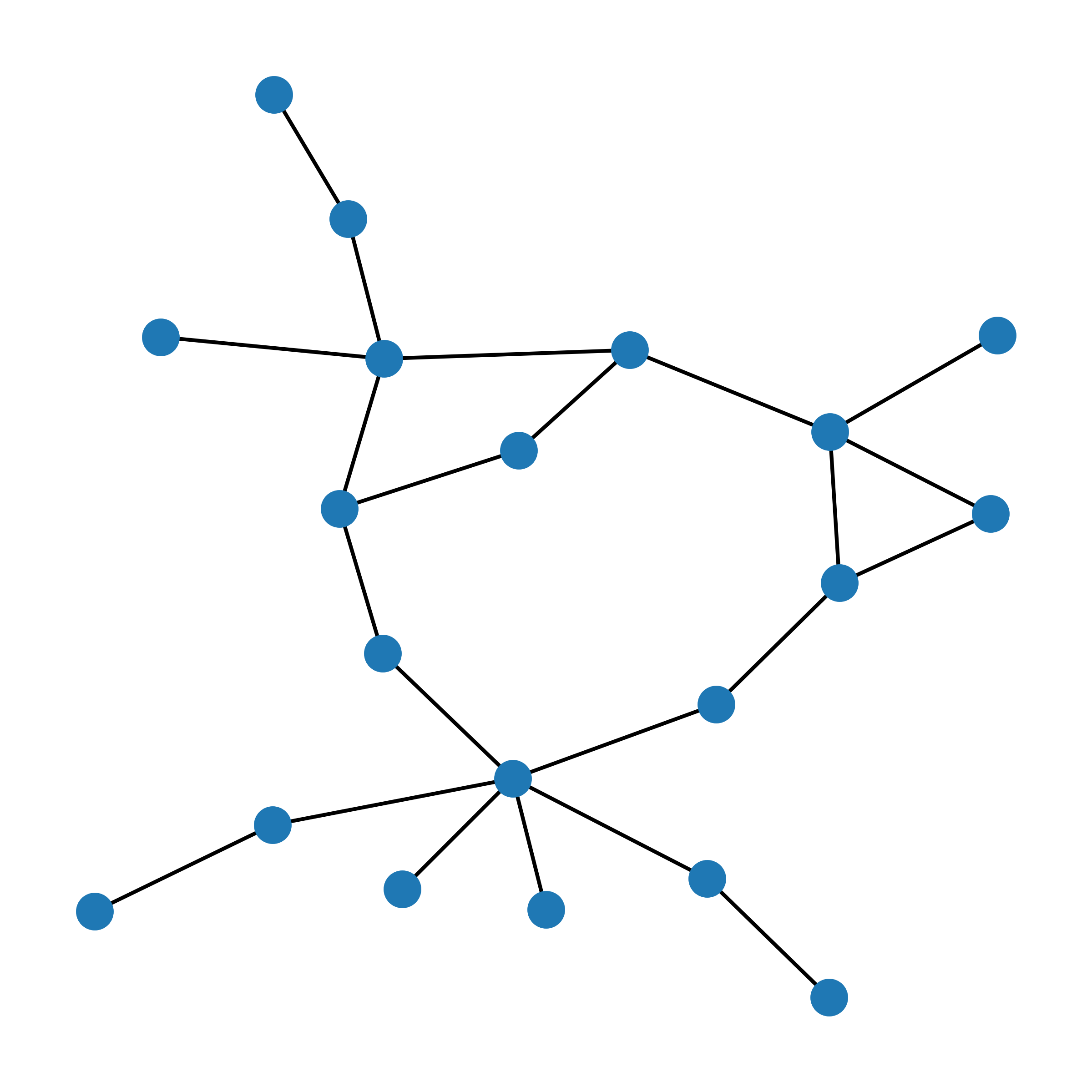}
  & \includegraphics[width=0.225\textwidth]{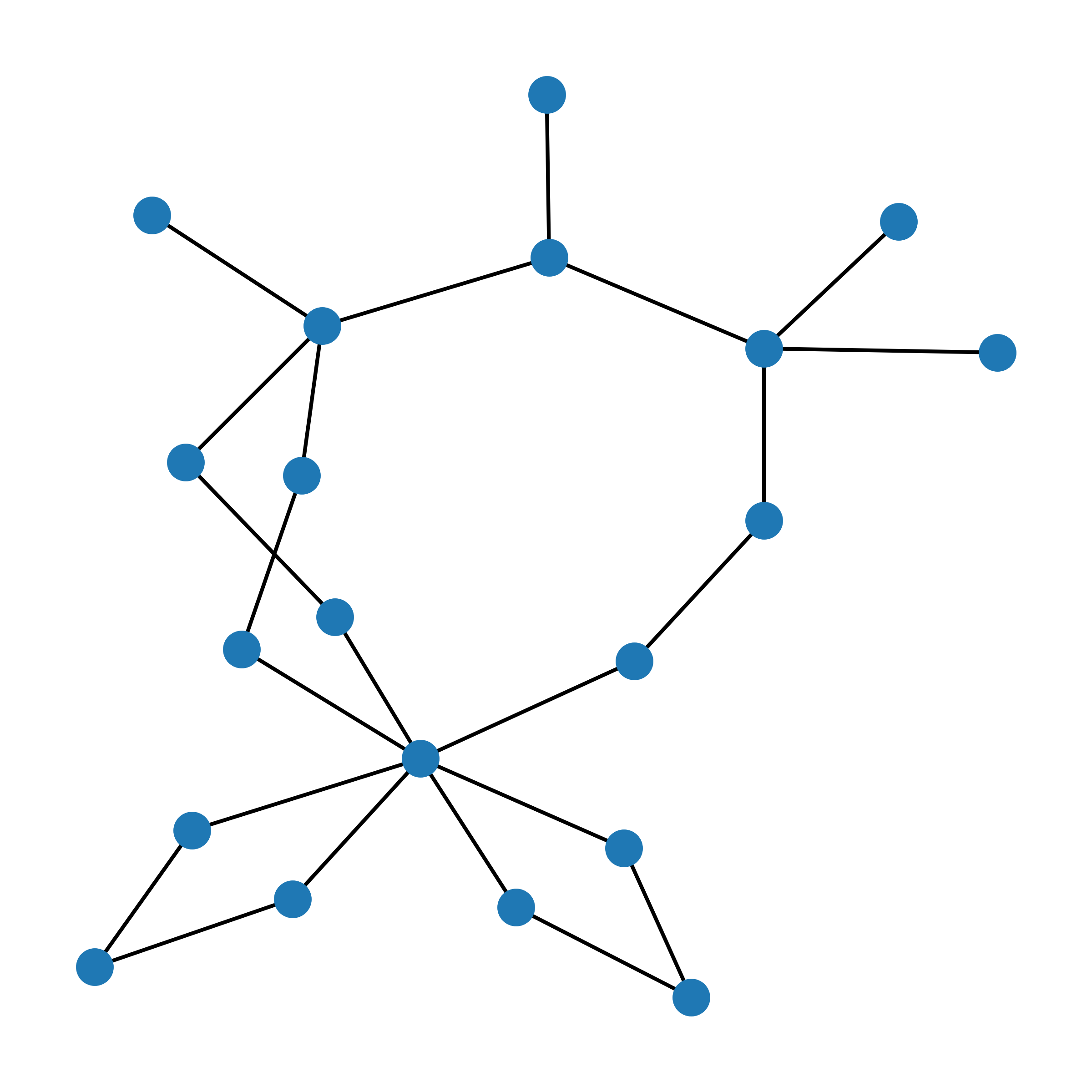} \\
gen=25, reward=-0.9167
  & gen=75, reward=-0.6717  
  & gen=275, reward=-0.1234 \\
\includegraphics[width=0.225\textwidth]{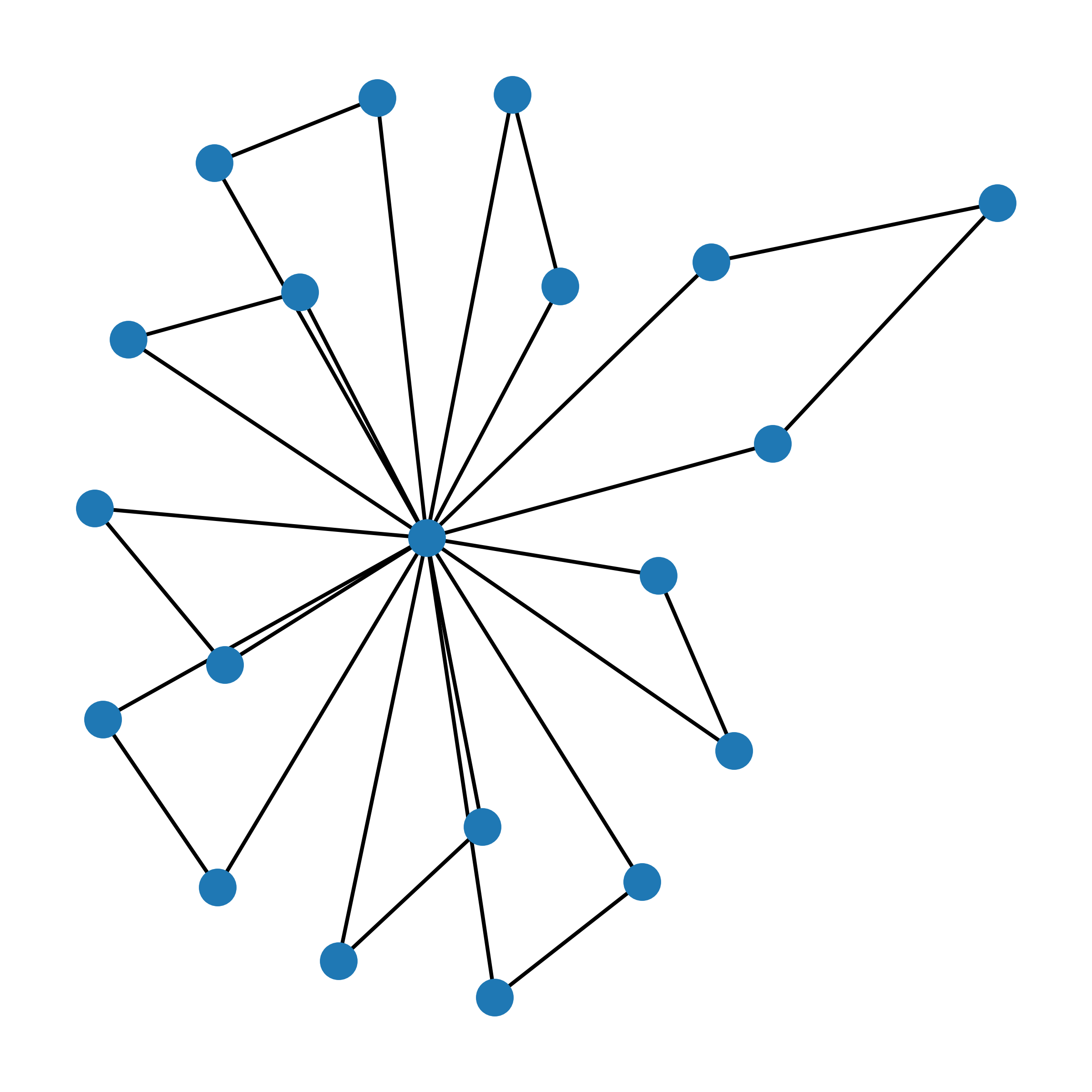}
  & \includegraphics[width=0.225\textwidth]{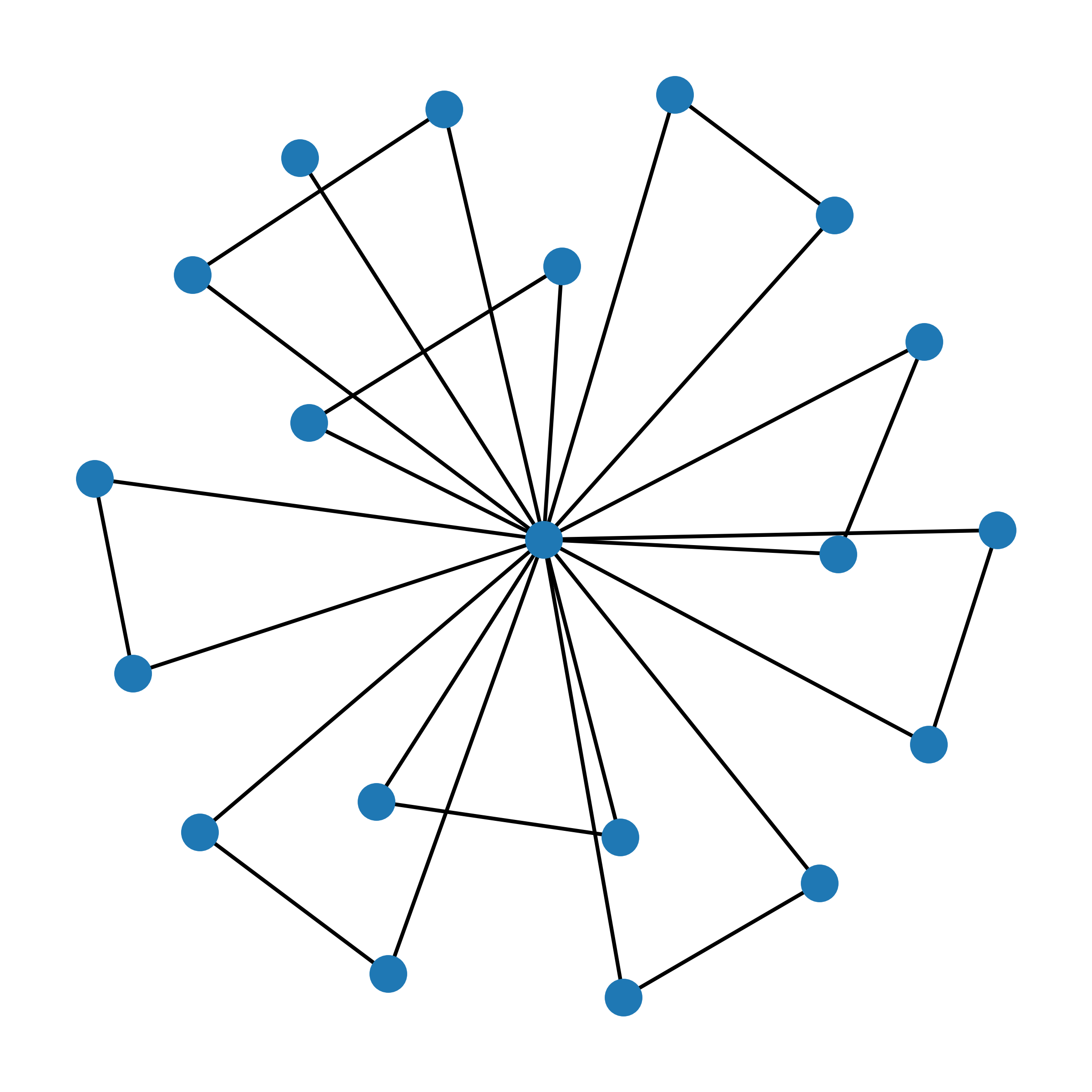} \\
gen=500, reward=2.2070
  & gen=925, reward=2.3460
\end{tabular}
\caption{Evolution of graphs with maximum rewards for the conjecture
         65: $\mu\leq\max_{v_i\sim v_j} \frac{(m_i+m_j)(d_im_i+d_jm_j)}{2m_im_j}$.}
\label{fig-evolution-65}
\end{figure}

Here we illustrate two particular conjectures
$$
31:\quad \mu\leq\max_{v\in V} \frac{4m_v^2}{m_v+d_v}
$$
and
$$
65:\quad \mu\leq\max_{v_i\sim v_j} \frac{(m_i+m_j)(d_im_i+d_jm_j)}{2m_im_j}
$$
for which RL managed to properly prove itself
by converging toward graphs with clearly identifiable and nontrivial structure.
Evolution of learning rewards for these two conjectures is shown in Fig.~\ref{fig-learning-evolution},
while evolutions of graphs with maximum rewards are shown in Figs. \ref{fig-evolution-31} and~\ref{fig-evolution-65}.
Learning rewards tend to oscillate after reaching a plateau
due to increases in action randomness,
which randomly adds or deletes edges to graphs constructed by the RL agent.
While the counterexamples for these two conjectures were found relatively early (after 250--500 generations),
we see that RL managed to further improve their rewards in later generations (after 600--900 generations).
A look at the maximum reward graphs shows that this prolonged learning was worthwile,
as it gave time to RL to make minor adjustments to its learning strategy and produce a more evident structure.

\begin{figure}
\centering
\begin{tabular}{ccc}
\includegraphics[width=0.225\textwidth]{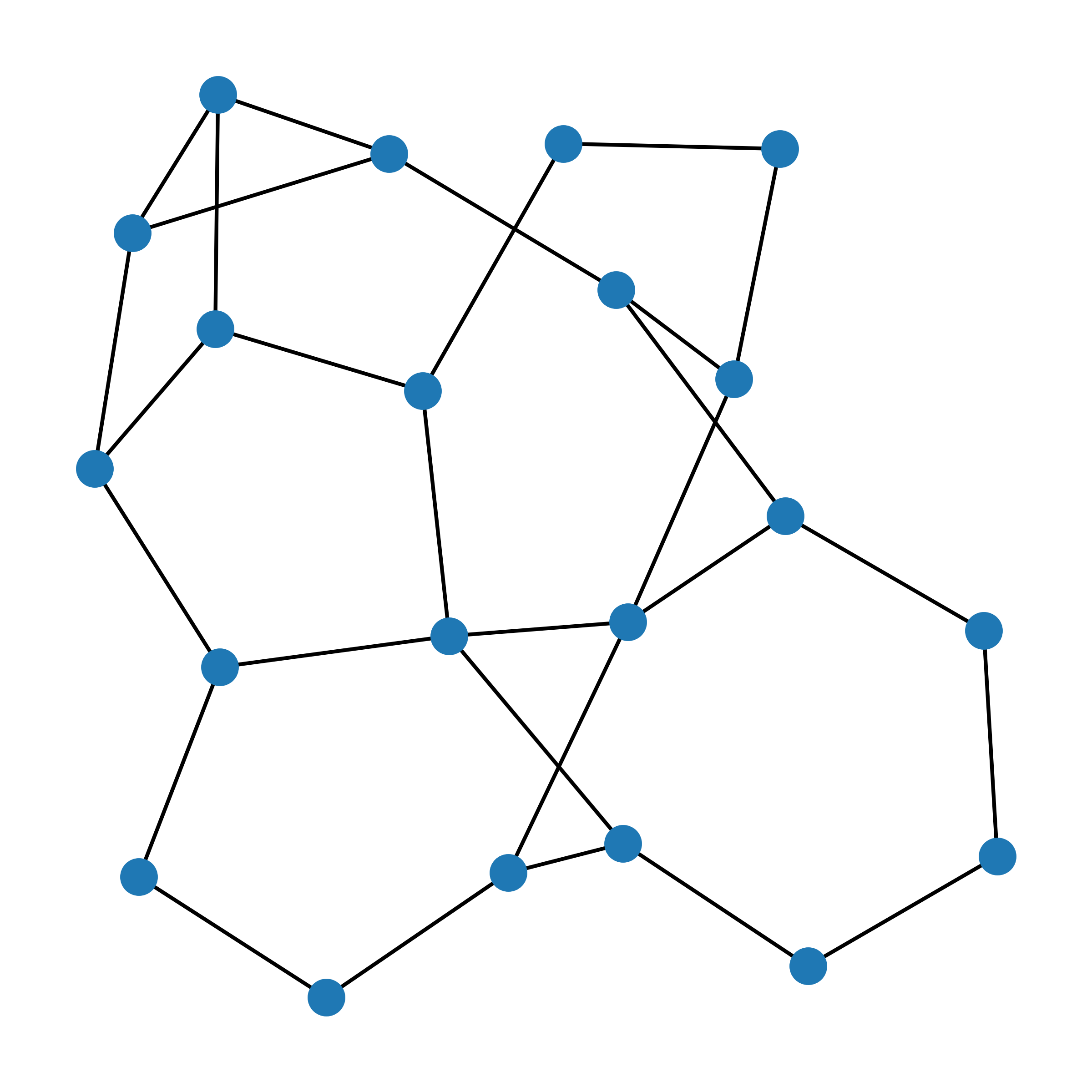} &
\includegraphics[width=0.225\textwidth]{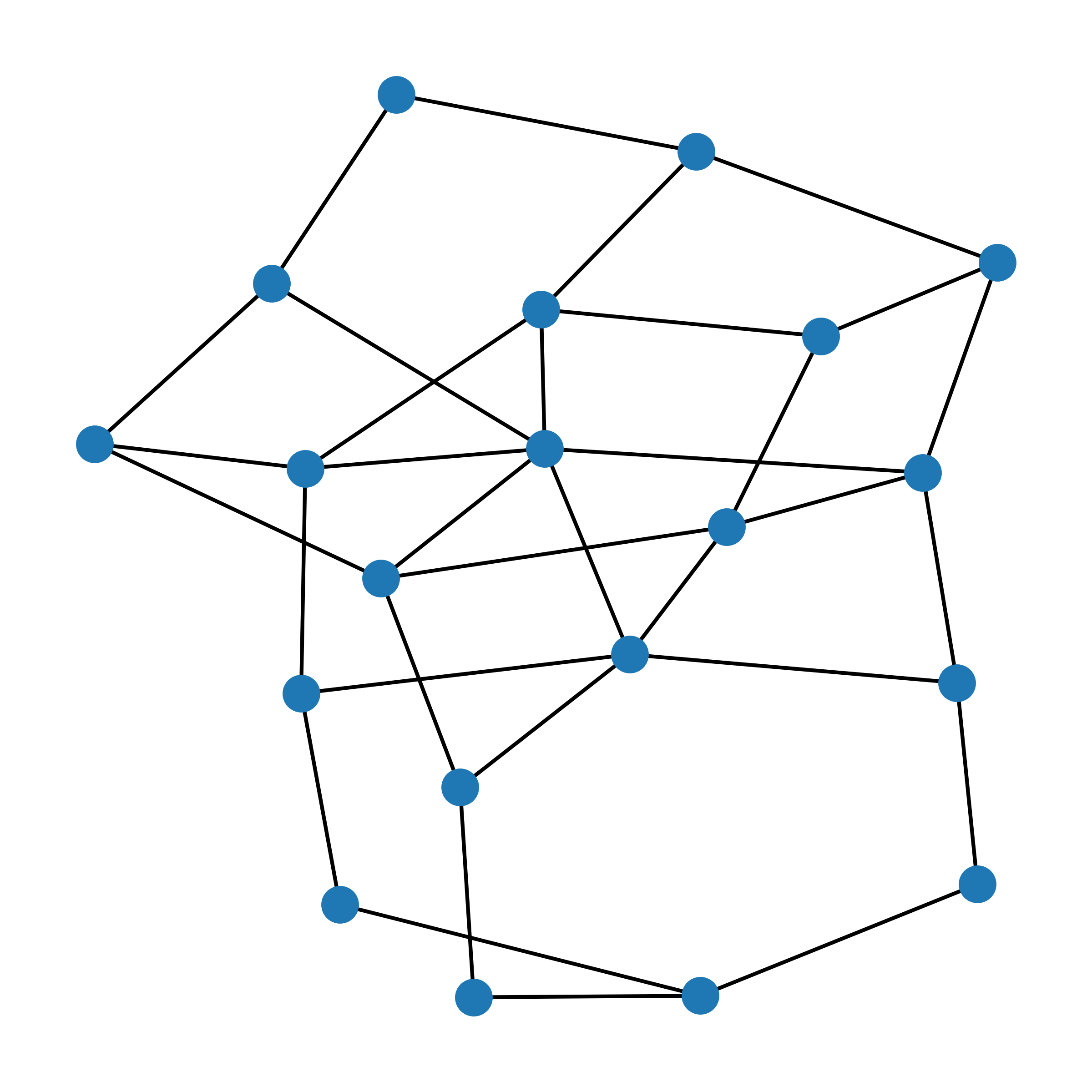} & 
\includegraphics[width=0.225\textwidth]{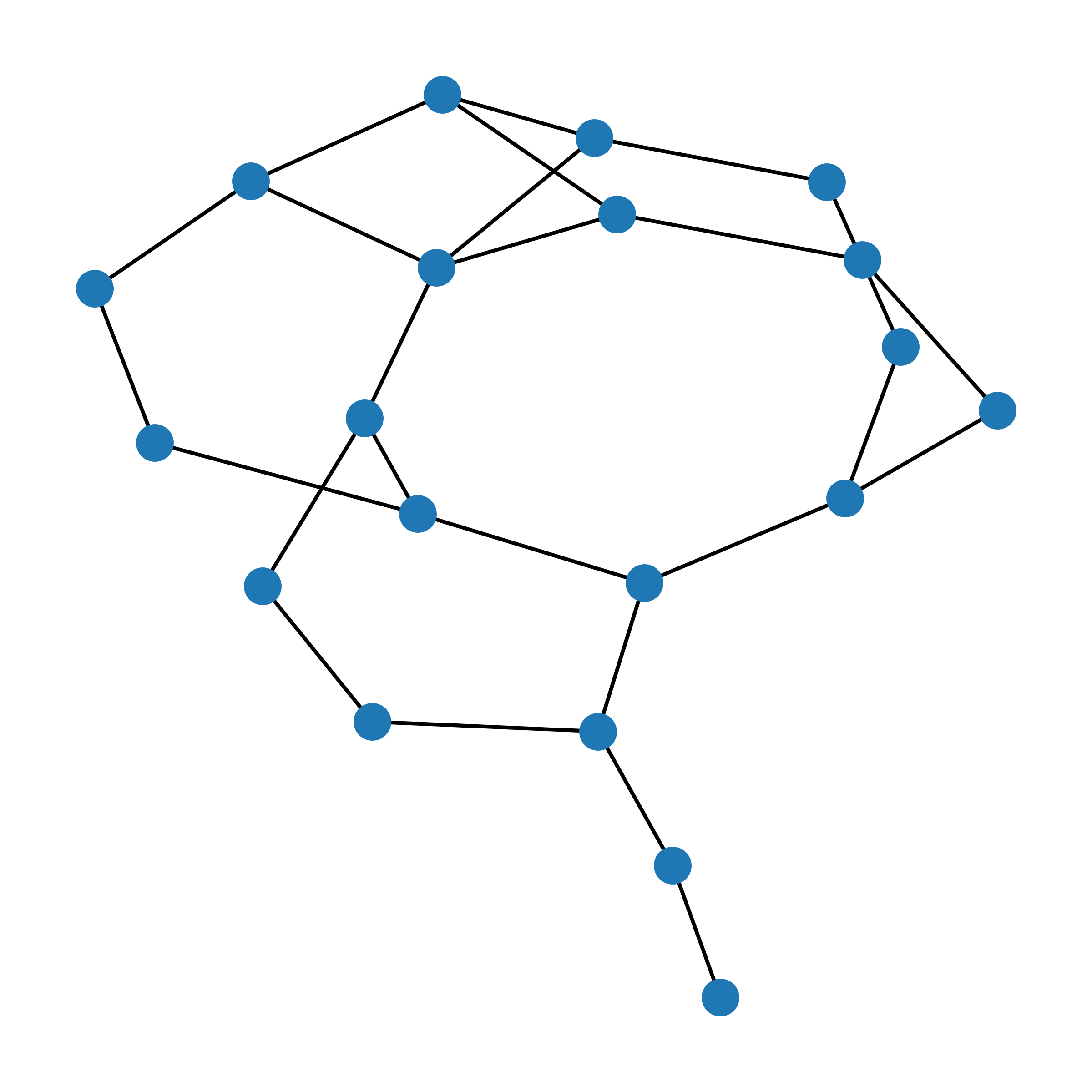} \\
Bound 3 &
Bound 15 & 
Bound 28 \\
\includegraphics[width=0.225\textwidth]{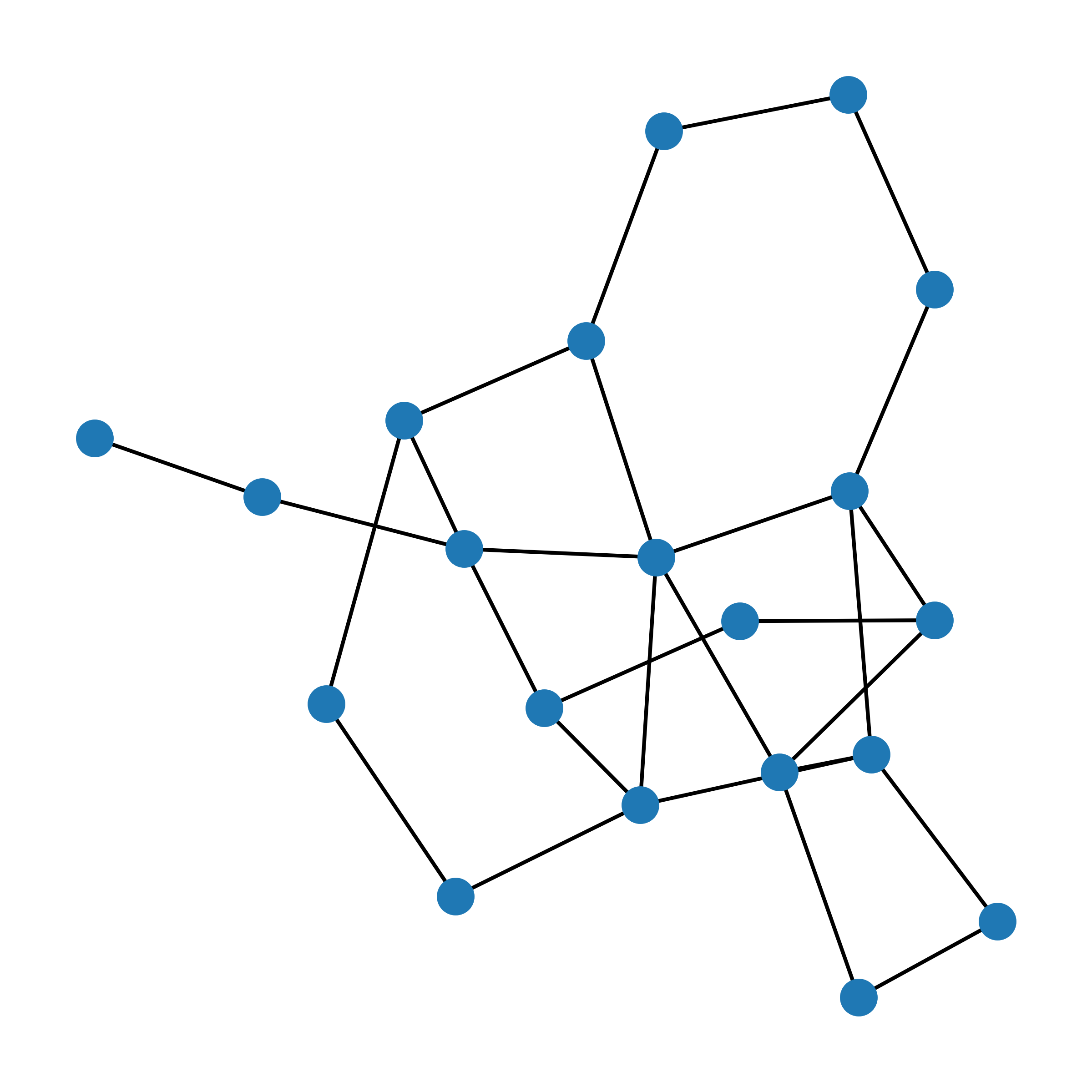} &
\includegraphics[width=0.225\textwidth]{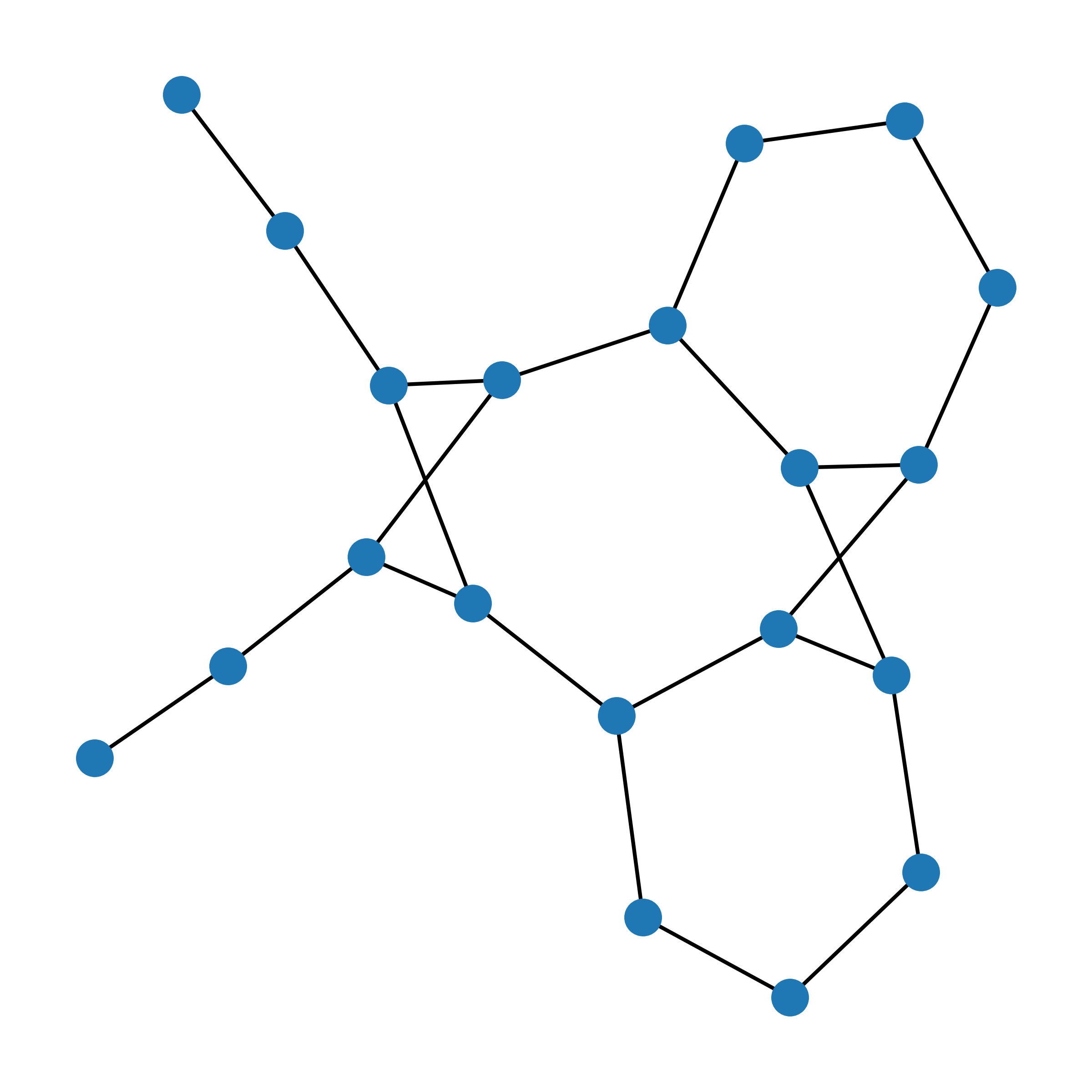} &
\includegraphics[width=0.225\textwidth]{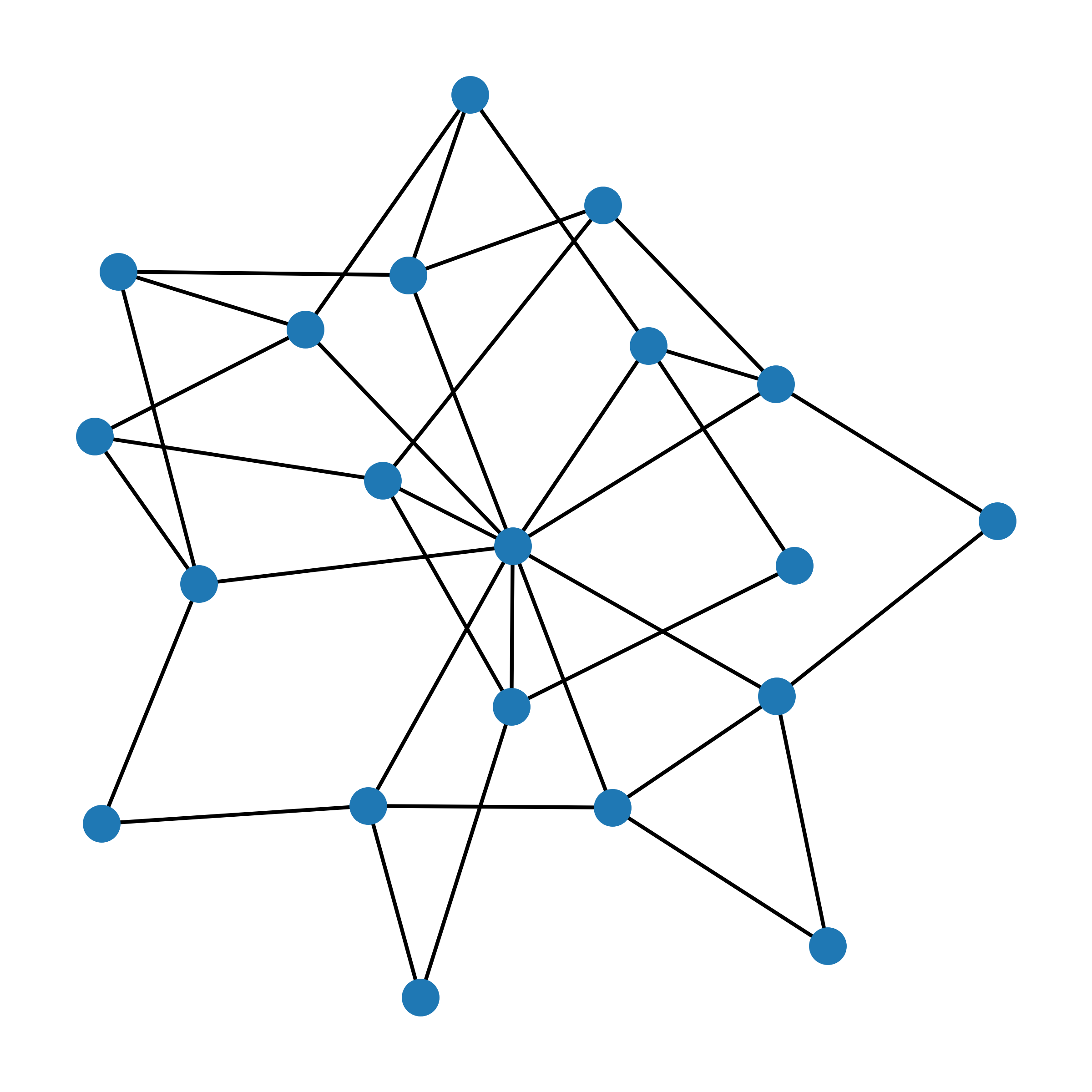} \\
Bound 29 &
Bound 36 &
Bound 49 \\
\includegraphics[width=0.225\textwidth]{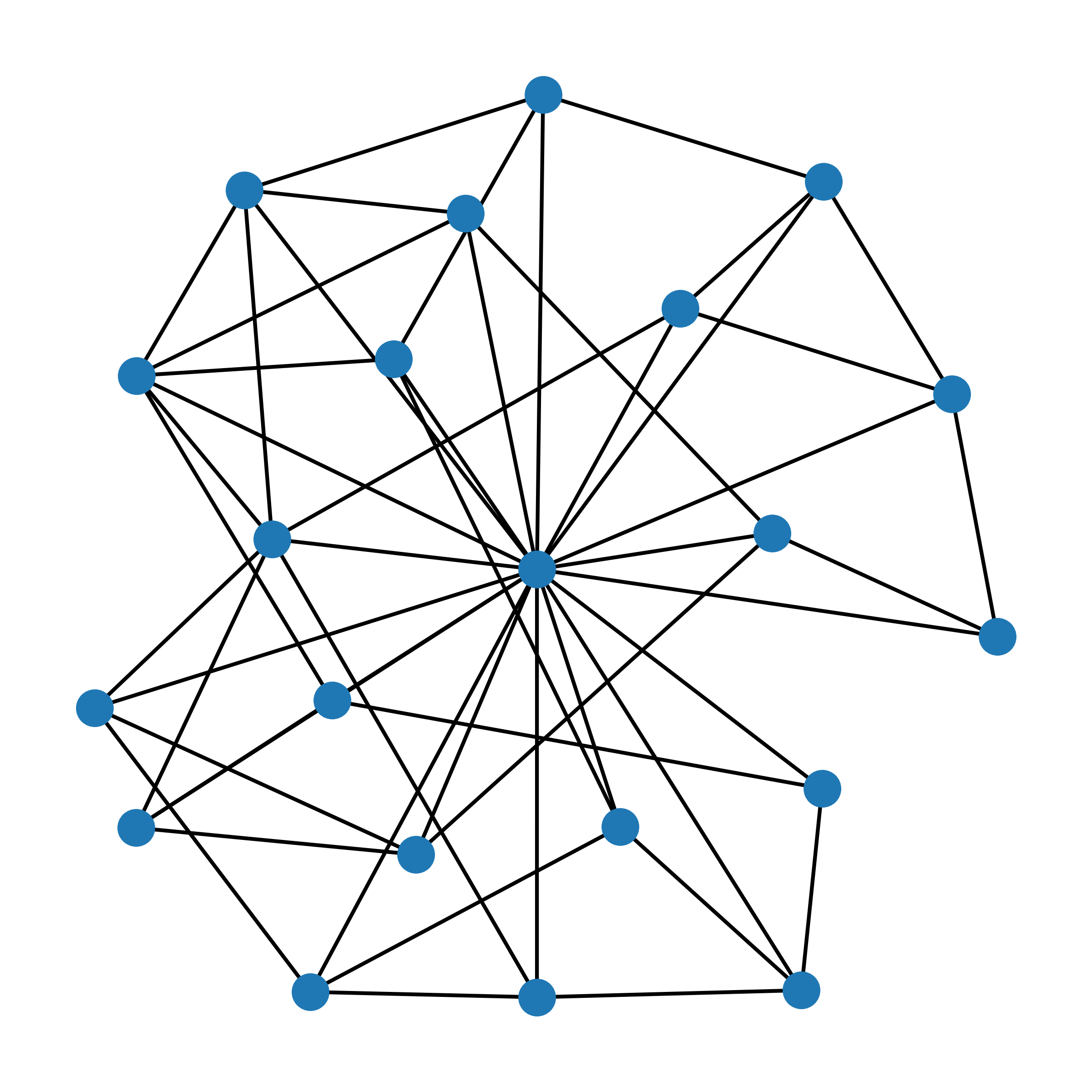} &
\includegraphics[width=0.225\textwidth]{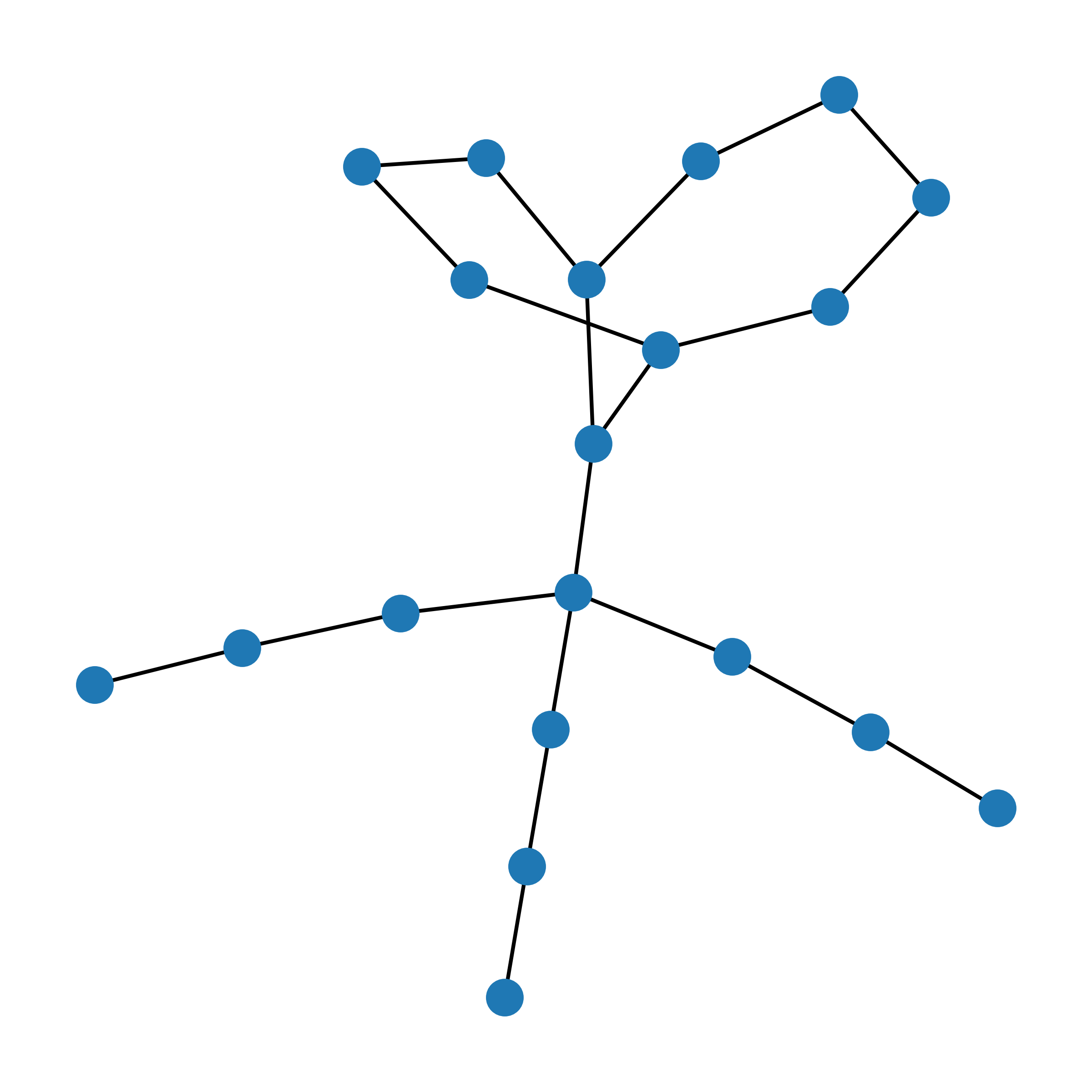} &
\includegraphics[width=0.225\textwidth]{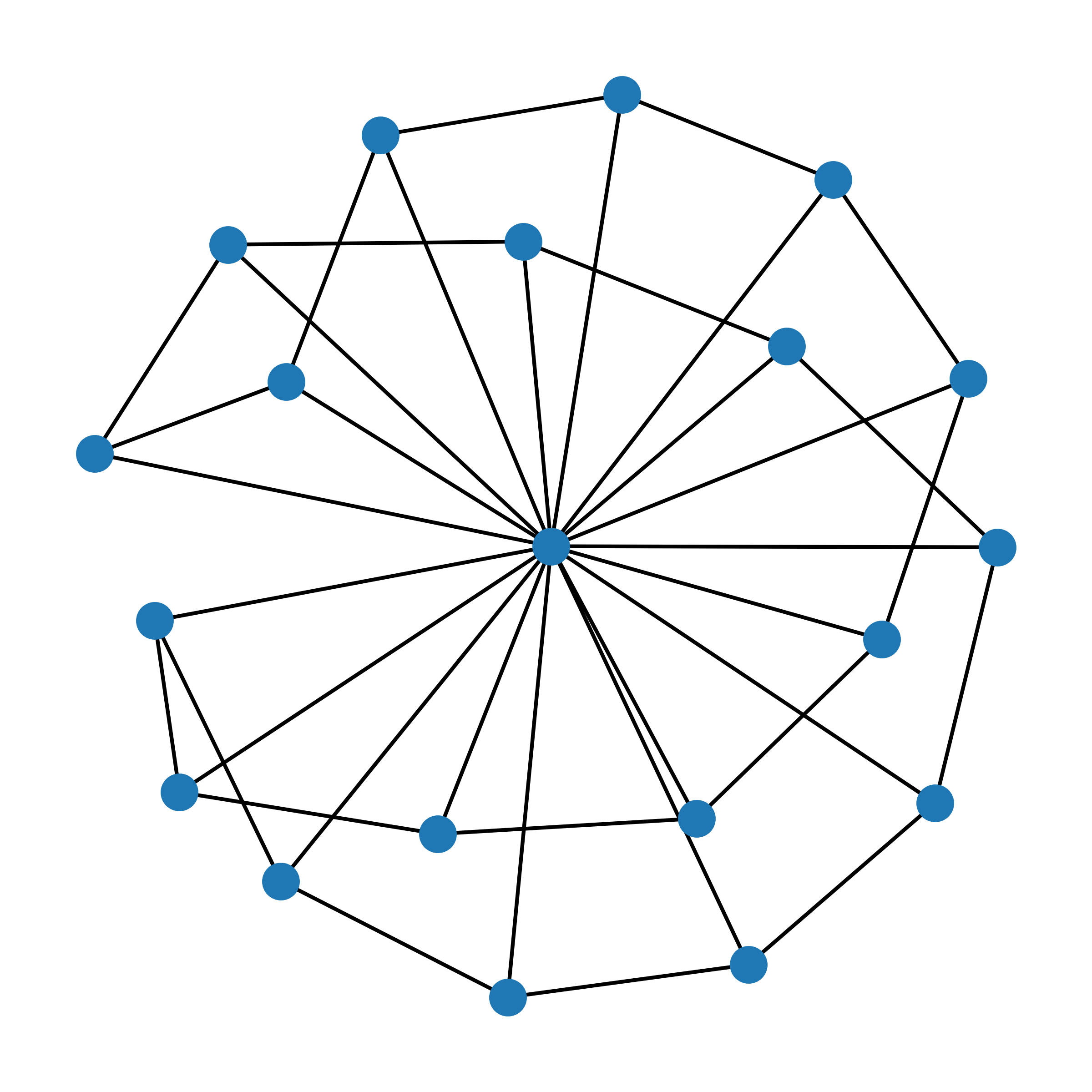} \\
Bound 51 &
Bound 67 &
Bound 68
\end{tabular}
\caption{Counterexamples for several other conjectured bounds from Appendix~\ref{sc-autolaplabounds-list}.}
\label{fig-other-counterexamples}
\end{figure}

Certainly, not all RL runs produce counterexamples with an easily identifiable structure,
as can be noted from the counterexamples shown in Fig.~\ref{fig-other-counterexamples}.
RL is, after all, an optimization algorithm that can get stuck in a local optimum.
Slow and steady increase of action randomness, which is incorporated in our reimplementation, 
may help it to get out of the local optimum if training is left to work over a larger number of generations.
While this was indeed the case with conjectured bounds 31 and~65,
there can be no upfront guarantee of success,
so it is up to the user to decide based on his/her intuition
when is the proper time to stop training that is unlikely to make further progress.

Nevertheless, the learning that does occur during the training of RL agent
can often suggest appropriate directions for further study.
For example,
while the counterexamples from Fig.~\ref{fig-other-counterexamples}
do not necessarily have clearly identifiable structure
(except obviously for a counterexample to bound 68),
they do suggest that counterexamples are most often subquartic graphs.

\begin{figure}
\centering
\begin{tabular}{cc}
\includegraphics[height=0.15\textheight]{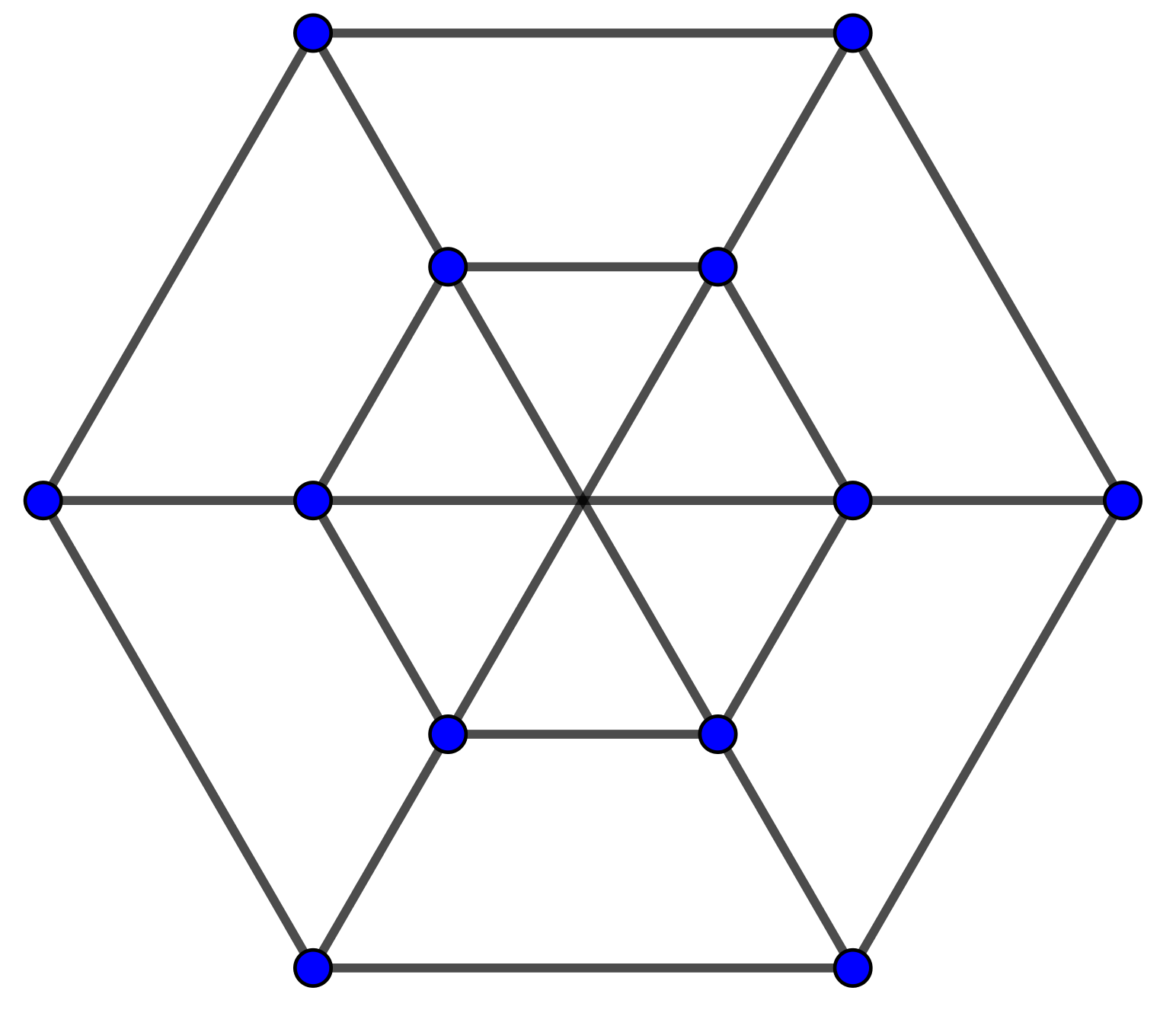} &
  \includegraphics[height=0.16\textheight]{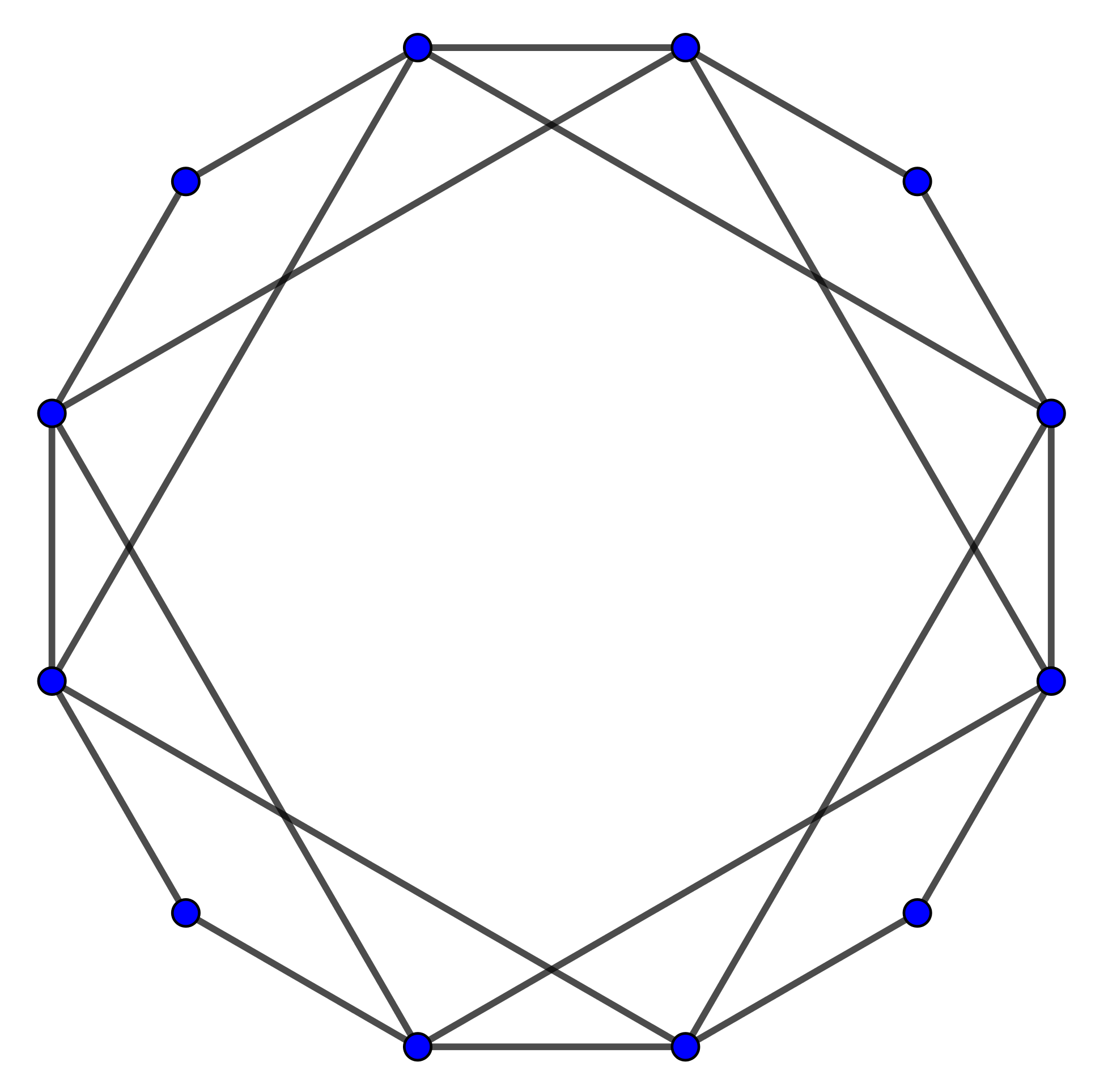} \\
$\mu(SQ^\ast)=7.41421$ &
  $\mu(SQ_{17})=7.37228$ \\
\\
\includegraphics[height=0.07\textheight]{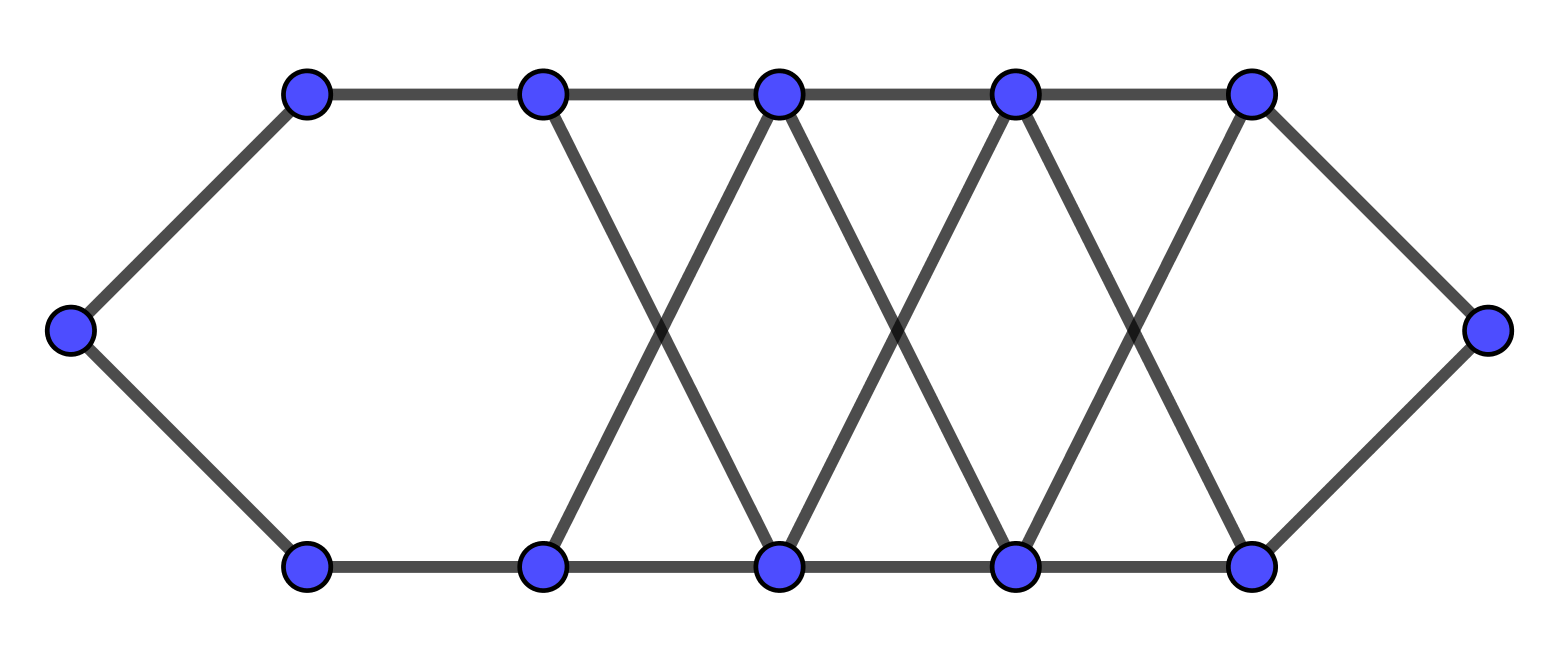} &
  \includegraphics[height=0.07\textheight]{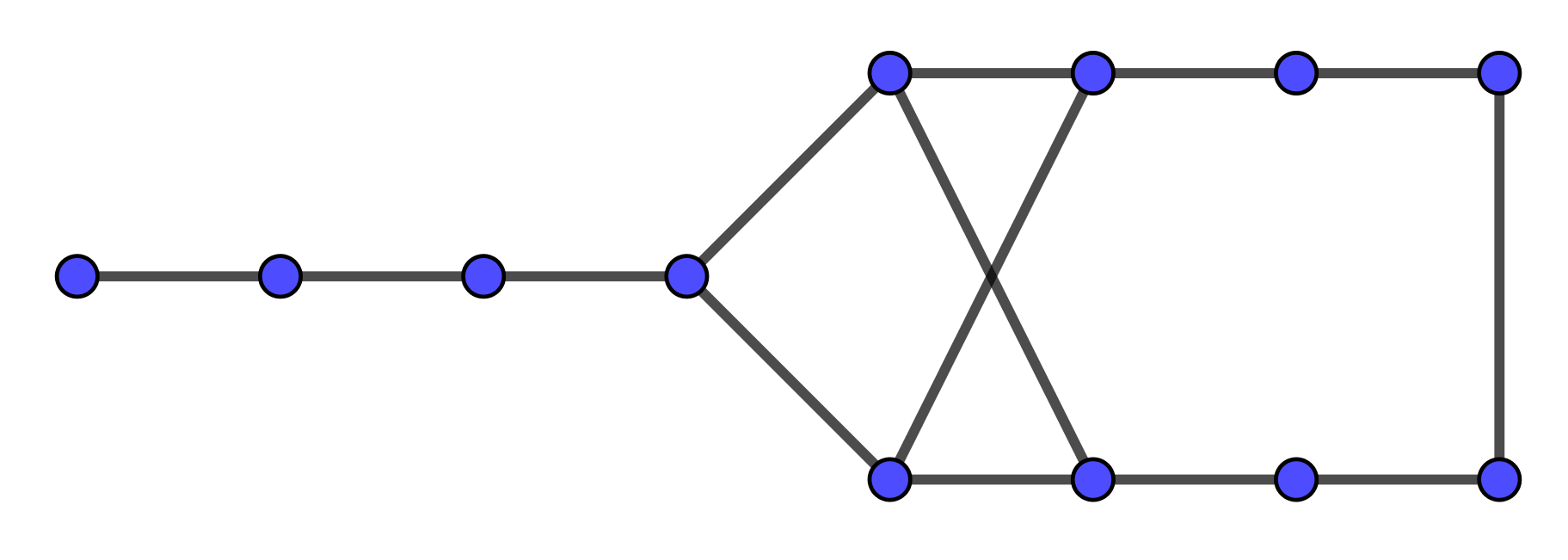} \\
$\mu(SQ_{50})=7.06459$ &
  $\mu(SQ_{66})=5.60286$
\end{tabular}
\caption{Counterexamples found among subquartic graphs on 12 vertices.
$SQ^\ast$ disproves 23 bounds: 2, 3, 15, 28, 29, 31, 32, 36, 43, 49, 52, 53, 54, 55, 57, 58, 59, 60, 61, 62, 63, 64 and 67,
while $SQ_{17}$, $SQ_{50}$ and $SQ_{66}$ disprove only the bound indicated by their subscripts.}
\label{fig-subquartic}
\end{figure}

Subquartic graphs are not too numerous, compared to connected graphs in general,
so we were able to perform an exhaustive search among subquartic graphs with up to 14~vertices.
It turned out that among them subquartic graphs on 12~vertices were most successful,
disproving 26 of the conjectured bounds,
including some for which RL did not find counterexamples,
as indicated in Table~\ref{tb-counterexamples}.
Interestingly, one of these graphs serves as a counterexample for a total of 23 of these bounds,
indicated as $SQ^\ast$ in Fig.~\ref{fig-subquartic}.
The remaining graphs in this figure represent counterexamples for bounds 17, 50 and~66.

On the other hand,
subquartic counterexamples have much more evident structure,
so that they are easily generalizable, 
together with the two RL counterexamples for bounds 31 and~65 
(see Figs. \ref{fig-evolution-31} and~\ref{fig-evolution-65}).
However, while these generalizations do yield several further examples of graphs
that disprove several conjectured bounds at once,
they do not manage to disprove any other bound than those already disproved by RL training and subquartic graphs,
so that we do not discuss them further here.

\definecolor{lg}{gray}{0.95}  
\rowcolors{3}{lg}{}  
\begin{table}[h!t]
\centering
\begin{tabular}{ccc}
\toprule
Conjectured & RL              & Subquartic      \\
bound       & counterexamples & counterexamples \\
\midrule
2  &            & $SQ^\ast$ \\
3  & \checkmark & $SQ^\ast$ \\
15 & \checkmark & $SQ^\ast$ \\
17 &            & $SQ_{17}$ \\
28 & \checkmark & $SQ^\ast$ \\
29 & \checkmark & $SQ^\ast$ \\
31 & \checkmark & $SQ^\ast$ \\
32 &            & $SQ^\ast$ \\
\midrule
36 & \checkmark & $SQ^\ast$ \\
41 & \checkmark &           \\
43 & \checkmark & $SQ^\ast$ \\
49 & \checkmark & $SQ^\ast$ \\
50 &            & $SQ_{50}$ \\
51 & \checkmark &           \\
52 & \checkmark & $SQ^\ast$ \\
53 & \checkmark & $SQ^\ast$ \\
54 & \checkmark & $SQ^\ast$ \\
55 & \checkmark & $SQ^\ast$ \\
57 & \checkmark & $SQ^\ast$ \\
58 & \checkmark & $SQ^\ast$ \\
59 & \checkmark & $SQ^\ast$ \\
60 & \checkmark & $SQ^\ast$ \\
61 &            & $SQ^\ast$ \\
62 & \checkmark & $SQ^\ast$ \\
63 & \checkmark & $SQ^\ast$ \\
64 & \checkmark & $SQ^\ast$ \\
65 & \checkmark &           \\
66 & \checkmark & $SQ_{66}$ \\
67 & \checkmark & $SQ^\ast$ \\
68 & \checkmark &           \\
\hiderowcolors  
\bottomrule
\end{tabular}
\caption{Disproved upper bounds on Laplacian spectral radius of graphs. 
(Refer to Appendix~\ref{sc-autolaplabounds-list} for the actual bounds indicated here.)}
\label{tb-counterexamples}
\end{table}

\section{Concluding remarks}

We have reimplemented here Wagner's approach of applying a cross-entropy method, 
a particular reinforcement learning technique,
to an environment that represents simple graphs,
in order to construct counterexamples in graph theory.

The reimplementation was applied to the set of 68 conjectured upper bounds 
on the Laplacian spectral radius of graphs from~\cite{BrHS},
which are listed in Appendix~\ref{sc-autolaplabounds-list}.
A total of 30 conjectured bounds was disproven:
25 directly by reinforcement learning 
and additional five by exhaustive search among subquartic graphs,
which was largely suggested by the properties of counterexamples obtained through reinforcement learning.
However, 38 conjectured bounds are still open after these computational attacks,
suggesting that their automated way of generation may have had some merit.
That appears to be especially true for vertex-based ones,
as only 8 out of 32 conjectured bounds were disproved here.
Probably due to the absence of symmetry in their expressions,
none of the newly conjectured bounds from~\cite{BrHS} has been proven in the last 18 years,
yet they still resist being disproved computationally.
We hope that new theoretical methods for proving at least a handful of these bounds will be discovered soon.

The reinforcement learning proved to be of valuable help in this study,
and it is definitely a promising area for further work,
either through the use of other available reinforcement learning methods,
or by developing new environments for more specific graph classes
or by applying it in new ways, e.g., to construct pairs of graphs at once in a single generation phase.

\appendix

\section{Invariants implemented in Graph class in graph6java}
\label{sc-graph6java}

The following list very briefly describes the invariants currently implemented in Graph class
of the Java archive graph6java.jar from 
\begin{quote}
\url{github.com/dragance106/cema-for-graphs}.
\end{quote}
Further details on implemented invariants can be found
in the source file of this class at \url{github.com/dragance106/graph6java}.

\begin{description}
\item[n()]   the number of vertices
\item[m()]   the number of edges

\item[degrees()]    the array of vertex degrees
\item[averageDegrees()]    the array of average degrees of neighbours of vertices

\item[isConnected()]    whether the graph is connected
\item[numberComponents()]    the number of connected components

\item[isIsomorphic(Graph h)]    whether the graph is isomorphic to another graph h
\item[matchingNumber()]    The matching number of the graph
\item[complement()]    the complement of the graph

\item[Amatrix()]    the adjacency matrix $A$ of the graph
\item[Acharpoly()]    the characteristic polynomial of $A$, 
                      returned as the array of $n+1$ coefficients from the highest at index 0 to the lowest at index $n$
\item[Aspectrum()]    the eigenvalues of $A$, sorted in non-decreasing order
\item[Aeigenvectors()]    the eigenvectors of corresponding eigenvalues of $A$, 
                          returned as a double array with eigenvectors placed in columns
\item[Acospectral(Graph h)]    whether the graph is cospectral to another graph h
\item[Aintegral()]    whether all eigenvalues of $A$ are integers
\item[Aenergy()]    the energy of $A$, 
                    understood as the absolute deviation of the eigenvalues of $A$ from their average value
                          
\item[Lmatrix()]    the Laplacian matrix
\item[Qmatrix()]    the signless Laplacian matrix
\item[Dmatrix()]    the distance matrix
\item[DLmatrix()]    the distance Laplacian matrix
\item[Mmatrix()]    the modularity matrix

Note that for each of the above matrices the Graph class also contains
the corresponding methods for computing its spectral properties,
i.e., Lcharpoly(), Qspectrum(), Deigenvectors(), DLenergy(), Mcospectral() etc.

\item[Asingular()]    whether 0 is an eigenvalue of $A$
\item[fiedlerVector()]    the eigenvector of the second smallest eigenvalue of $L$
\item[LEL()]   Laplacian-like-energy invariant
\item[estrada()]    Estrada index
\item[Lestrada()]    Laplacian Estrada index

\item[diameter()]    the diameter
\item[radius()]   the radius
\item[wiener()]   Wiener index
\item[transmissions()]    the array of vertex transmissions, i.e.,
                          the sums of distances from a vertex to all other vertices
\item[transmissionIrregular()]    whether the graph is transmission irregular
\item[szeged()]    Szeged index
\item[weightedSzeged()]   weighted Szeged index
\item[randic()]    Randi\'c index
\item[zagreb1()]    The first Zagreb index
\item[zagreb2()]    The second Zagreb index
\end{description}

\section{The list of conjectured upper bounds from~\cite{BrHS}}
\label{sc-autolaplabounds-list}

The first group of bounds is of the form
$$
\mu\leq\max_{v\in V} f(d_v, m_v)
$$
where $\mu$ is the spectral radius of Laplacian matrix of graph~$G$,
$V$ is its vertex set,
while $d_v$ and $m_v$ are, respectively, the degree of~$v$ and the average degree of the neighbours of~$v$, for $v\in V$.
In the following list we present only the right-hand side of these bounds.
The symbol O after the ordinal number means that the conjecture is still open,
while the symbol X means that the counterexample has been found in this paper,
as indicated in Table~\ref{tb-counterexamples}.

{\small
\begin{longtable}{rlrl}
1. O & $\max_{v\in V}\sqrt{\frac{4 d_v^3}{m_v}}$
  & 17. X & $\max_{v\in V}\sqrt[4]{5 d_v^4 + 11 m_v^4}$ \\
2. X & $\max_{v\in V}\frac{2 m_v^2}{d_v}$
  & 18. O & $\max_{v\in V}\sqrt{\frac{2 m_v^3}{d_v} + 2 d_v^2}$ \\
3. X & $\max_{v\in V}\frac{m_v^2}{d_v} + m_v$
  & 19. O & $\max_{v\in V}\sqrt[4]{4 d_v^4 + 12 d_vm_v^3}$ \\
4. O & $\max_{v\in V}\frac{2 d_v^2}{m_v}$
  & 20. O & $\max_{v\in V}\frac{\sqrt{7 d_v^2 + 9 m_v^2}}{2}$ \\
5. O & $\max_{v\in V}\frac{d_v^2}{m_v} + m_v$
  & 21. O & $\max_{v\in V}\sqrt{\frac{d_v^3}{m_v} + 3 m_v^2}$ \\
6. O & $\max_{v\in V}\sqrt{m_v^2 + 3 d_v^2}$
  & 22. O & $\max_{v\in V}\sqrt[4]{2 d_v^4 + 14 d_v^2m_v^2}$ \\
7. O & $\max_{v\in V}\frac{d_v^2}{m_v} + d_v$
  & 23. O & $\max_{v\in V}\sqrt{d_v^2 + 3 d_vm_v}$ \\
8. O & $\max_{v\in V}\sqrt{d_v(m_v + 3 d_v)}$
  & 24. O & $\max_{v\in V}\sqrt[4]{6 d_v^4 + 10 m_v^4}$ \\
9. O & $\max_{v\in V}\frac{m_v + 3 d_v}{2}$
  & 25. O & $\max_{v\in V}\sqrt[4]{3 d_v^4 + 13 d_v^2m_v^2}$ \\
10. O & $\max_{v\in V}\sqrt{d_v(d_v + 3 m_v}$
  & 26. O & $\max_{v\in V}\frac{\sqrt{5 d_v^2 + 11 d_vm_v}}{2}$ \\
11. O & $\max_{v\in V}\frac{2 m_v^3}{d^{2}_v}$
  & 27. O & $\max_{v\in V}\sqrt{\frac{3 d_v^2 + 5 d_vm_v}{2}}$ \\
12. O & $\max_{v\in V}\sqrt{2 m_v^2 + 2 d_v^2}$
  & 28. X & $\max_{v\in V}\sqrt{\frac{2 m_v^4}{d_v^2} + 2 d_vm_v}$ \\
13. O & $\max_{v\in V}\frac{2 m_v^4}{d_v^3}$
  & 29. X & $\max_{v\in V}\sqrt{m_v^2 + \frac{3 m_v^3}{d_v}}$ \\
14. O & $\max_{v\in V}\frac{2 d_v^3}{m_v^2}$
  & 30. O & $\max_{v\in V}\frac{m_{v}^{3}}{d_{v}^{2}}+\frac{d_v^2}{m_v}$ \\
15. X & $\max_{v\in V}\sqrt{\frac{4 m_v^3}{d_v}}$
  & 31. X & $\max_{v\in V}\frac{4 m_v^2}{m_v + d_v}$ \\
16. O & $\max_{v\in V}\frac{2 d_v^4}{m_v^3}$
  & 32. X & $\max_{v\in V}\frac{\sqrt{m_v^3 (m_v + 3 d_v)}}{d_v}$
\end{longtable}
}

The second group of bounds is of the form
$$
\mu\leq\max_{v_i\sim v_j} f(d_{v_i}, m_{v_i}, d_{v_j}, m_{v_j})
$$
where the maximum is taken over all pairs of adjacent vertices $v_i$, $v_j$ in~$V$.
Again, in the following list we present only the right-hand side of these bounds,
with $d_{v_i}, m_{v_i}, d_{v_j}, m_{v_j}$ shortened as $d_i, m_i, d_j, m_j$, respectively.
{\small
\begin{longtable}{rl}
33. O & $\max_{v_i\sim v_j}2(d_i + d_j) -(m_i + m_j)$ \\
34. O & $\max_{v_i\sim v_j}\frac{2(d_i^2 + d_j^2)}{d_i + d_j}$ \\
35. O & $\max_{v_i\sim v_j}\frac{2(d_i^2 + d_j^2)}{m_i + m_j}$ \\
36. X & $\max_{v_i\sim v_j}\frac{2(m_i^2 + m_j^2)}{d_i + d_j}$ \\
37. O & $\max_{v_i\sim v_j}\sqrt{2(d_i^2 + d_j^2)}$ \\
38. O & $\max_{v_i\sim v_j}2 + \sqrt{2(d_i-1)^2 + 2(d_j-1)^2}$ \\
39. O & $\max_{v_i\sim v_j}2 + \sqrt{2(d_i^2 + d_j^2) -4(m_i + m_j) + 4}$ \\
40. O & $\max_{v_i\sim v_j}2 + \sqrt{2\left((m_i-1)^2 + (m_j-1)^2\right) + (d_i^2 + d_j^2) -(d_im_i + d_jm_j)}$ \\
41. X & $\max_{v_i\sim v_j}2 + (m_i + m_j) -(d_i + d_j) +\sqrt{2(d_i^2 + d_j^2) -4(m_i + m_j) + 4}$ \\
42. O & $\max_{v_i\sim v_j}\sqrt{d_i^2 + d_j^2 + 2 m_im_j}$ \\
43. X & $\max_{v_i\sim v_j}2 + \sqrt{3(m_i^2 + m_j^2) -2 m_im_j -4(d_i + d_j) + 4}$ \\
44. O & $\max_{v_i\sim v_j}2 + \sqrt{2\left((d_i-1)^2 + (d_j-1)^2 + m_im_j - d_id_j\right)}$ \\
45. O & $\max_{v_i\sim v_j}2 + \sqrt{(d_i - d_j)^2 + 2(d_im_i + d_jm_j) -4(m_i + m_j) + 4}$ \\
46. O & $\max_{v_i\sim v_j}2 + \sqrt{2(d_i^2 + d_j^2) -16 \frac{d_id_j}{m_i + m_j} + 4}$ \\
47. O & $\max_{v_i\sim v_j}\frac{2(d_i^2 + d_j^2) -(m_i-m_j)^2}{d_i + d_j}$ \\
48. O & $\max_{v_i\sim v_j}\frac{2(d_i^2 + d_j^2)}{2 + \sqrt{2(d_i^2 + d_j^2) -4(m_i + m_j) + 4}}$ \\
49. X & $\max_{v_i\sim v_j}2 + \sqrt{2(m_i^2 + m_j^2) + (d_i-d_j)^2 -4(d_i + d_j) + 4}$ \\
50. X & $\max_{v_i\sim v_j}2\frac{d_i^2 + d_j^2 + m_im_j - d_id_j}{d_i + d_j}$ \\
51. X & $\max_{v_i\sim v_j}2(m_i + m_j) -4 \frac{m_im_j}{d_i + d_j}$ \\
52. X & $\max_{v_i\sim v_j}2 + \sqrt{\sqrt{8(m_i^4 + m_j^4) -8(d_i^2 + d_j^2) + 4} -4(d_i + d_j) + 6}$ \\
53. X & $\max_{v_i\sim v_j}2 + \sqrt{\sqrt{8(m_i^4 + m_j^4) -8(d_im_i + d_jm_j) + 4} -4(d_i + d_j) + 6}$ \\
54. X & $\max_{v_i\sim v_j}2 + \sqrt{2(m_i^2 + m_j^2) + (d_im_i + d_jm_j) -(d_i^2 + d_j^2) -4(d_i + d_j) + 4}$ \\
55. X & $\max_{v_i\sim v_j}2 + \sqrt{3(m_i^2 + m_j^2) -(d_i^2 + d_j^2) -4(m_i + m_j) + 4}$ \\
56. O & $\max_{v_i\sim v_j}\frac{(d_i^2 + d_j^2)(m_i + m_j)}{2d_id_j}$ \\
57. X & $\max_{v_i\sim v_j}2 + \sqrt{2(m_i^2 + m_j^2) -8 \frac{d_i^2 + d_j^2}{m_i + m_j} + 4}$ \\
58. X & $\max_{v_i\sim v_j}2 + \sqrt{2(m_i^2 + m_im_j + m_j^2) -(d_im_i + d_jm_j) -4(d_i + d_j) + 4}$ \\
59. X & $\max_{v_i\sim v_j}\frac{2(m_i^2 + m_im_j + m_j^2) -(d_i^2 + d_j^2)}{m_i + m_j}$ \\
60. X & $\max_{v_i\sim v_j}2 + \sqrt{2(m_i^2 + m_im_j + m_j^2) -(d_i^2 + d_j^2) -4(d_i + d_j) + 4}$ \\
61. X & $\max_{v_i\sim v_j}\frac{2(m_i^2 + m_j^2)}{2 + \sqrt{2\left((d_i-1)^2 + (d_j-1)^2\right)}}$ \\
62. X & $\max_{v_i\sim v_j}2 + \sqrt{m_i^2+ 4 m_im_j + m_j^2 -2 d_id_j -4(d_i + d_j) + 4}$ \\
63. X & $\max_{v_i\sim v_j}d_i + d_j + m_i + m_j -4\frac{d_id_j}{m_i + m_j}$ \\
64. X & $\max_{v_i\sim v_j}\frac{m_im_j(d_i + d_j)}{d_id_j}$ \\
65. X & $\max_{v_i\sim v_j}\frac{(m_i + m_j)(d_im_i + d_jm_j)}{2m_im_j}$ \\
66. X & $\max_{v_i\sim v_j}\frac{m_i^2 + 4 m_im_j + m_j^2 -(d_im_i + d_jm_j)}{d_i + d_j}$ \\
67. X & $\max_{v_i\sim v_j}\frac{(m_i + m_j)(d_im_i + d_jm_j)}{2d_id_j}$ \\
68. X & $\max_{v_i\sim v_j}2 + \sqrt{(m_i-m_j)^2 + 4 d_id_j -4(m_i + m_j) + 4}$
\end{longtable}
}

\end{document}